\documentclass[11pt,reqno,oneside]{amsart}
\usepackage[T1]{fontenc}
\usepackage{lmodern,amssymb,mathtools,microtype}
\usepackage[a4paper,top=25mm,bottom=25mm,left=25mm,right=25mm]{geometry}
\usepackage[hidelinks,bookmarksnumbered=true]{hyperref}
\numberwithin{equation}{section}
\newtheorem{theorem}{Theorem}[section]
\newtheorem{lemma}[theorem]{Lemma}
\newtheorem{proposition}[theorem]{Proposition}

\theoremstyle{remark}

\DeclareMathOperator{\tr}{tr}
\DeclareMathOperator{\diver}{div}

\newcommand{\E}{\mathbb E}
\newcommand{\Sph}{\mathbb S}
\newcommand{\doi}[1]{\href{https://doi.org/#1}{\nolinkurl{doi:#1}}}
\newcommand{\arxiv}[1]{\href{https://arxiv.org/abs/#1}{\nolinkurl{arXiv:#1}}}
\makeatletter
\renewcommand{\subsection}{\@startsection{subsection}{2}{\z@}%
  {12pt plus 3pt minus 2pt}{5pt plus 1pt}{\normalfont\bfseries}}
\renewcommand{\paragraph}{\@startsection{paragraph}{4}{\z@}%
  {8pt plus 2pt minus 1pt}{-.5em}{\normalfont\itshape}}
\renewcommand{\@setauthors}{%
  \begin{center}
    \vspace{1em}%
    \normalfont\small
    \authors\\
    \authoraffiliation
  \end{center}}
\makeatother
\allowdisplaybreaks[2]
\hypersetup{pdftitle={Yau's first eigenvalue conjecture for embedded minimal hypersurfaces with constant scalar curvature},pdfauthor={Ruihan Chen},pdfkeywords={Yau conjecture, minimal hypersurface, first eigenvalue, constant scalar curvature}}
\title[Yau's conjecture for hypersurfaces with constant scalar curvature]{Yau's first eigenvalue conjecture for embedded minimal hypersurfaces with constant scalar curvature}
\author{Ruihan Chen}
\newcommand{\authoraffiliation}{Faculty of Science, The University of Hong Kong}
\date{}
\subjclass[2020]{Primary 53C42, 58J50. Secondary 35P15, 58J35}
\keywords{Minimal hypersurface, Yau's first eigenvalue conjecture, constant scalar curvature, single-layer potential, spherical heat kernel}
\begin{document}
\begin{abstract}
Let $M^n$ be a smooth closed connected embedded minimal hypersurface of the unit sphere.
Suppose that $M$ has constant scalar curvature. We prove that its first positive
Laplace eigenvalue is $n$ and that the ambient coordinate functions span its
first eigenspace. To prove both assertions, we consider an eigenfunction $f$
with eigenvalue $0<\lambda\le n$ orthogonal to every coordinate function.
We use $f$ as the density of a single-layer potential on the sphere.
The geometry of $M$ gives a lower bound for the Dirichlet energy of this
potential. The spherical heat kernel gives a smaller upper bound for the same
energy. This contradiction rules out every such nonzero $f$.
\end{abstract}
\maketitle

\section{Introduction}
Let $M^n$ be a smooth closed connected embedded minimal hypersurface of the
unit sphere, and let
\[
  X\colon M^n\hookrightarrow\Sph^{n+1}(1)\subset\mathbb R^{n+2}
\]
denote the inclusion. Write $A$ for the
shape operator of $M$ and $S=|A|^2$ for the sum of the squares of its principal
curvatures. We use the nonnegative Laplacian $L=-\Delta_M$, where
$\Delta_M=\diver\nabla$. We list its eigenvalues with multiplicity as
\[
  0=\lambda_0<\lambda_1\le\lambda_2\le\cdots.
\]

Write $X=(X^1,\ldots,X^{n+2})$ for the position vector. By Takahashi's
theorem, minimality gives $LX=nX$ \cite{Takahashi}. Thus each nonzero
coordinate function $X^i$ is an eigenfunction with eigenvalue $n$.
Yau's first eigenvalue conjecture asserts that $\lambda_1=n$ for every such hypersurface
\cite[Problem~100]{Yau}, \cite[p.~562]{ChoiWang}. Kusner and Wang asked
whether the coordinate functions span the whole first eigenspace
\cite[Conjecture~4]{KM}. We consider both questions when $S$ is constant.

\begin{theorem}\label{main}
Let $M^n\subset\Sph^{n+1}(1)$ be a smooth closed connected embedded
minimal hypersurface, $n\ge1$. If $S=|A|^2$ is constant, then
\[
  \lambda_1=n,\qquad
  \ker(L-n)=V_X,\qquad V_X=\operatorname{span}\{X^1,\ldots,X^{n+2}\}.
\]
If $M$ is an equator, then $\dim V_X=n+1$. If $M$ is not an equator,
then $\dim V_X=n+2$.
\end{theorem}

The Gauss equation gives
\[
  \operatorname{Scal}_M=n(n-1)-S.
\]
Thus $S$ is constant exactly when $M$ has constant scalar curvature. Lei, Xu and Xu
asked whether this condition implies $\lambda_1=n$
\cite[Problem~4.1(i)]{LXX}. Theorem~\ref{main} answers their question
for embedded hypersurfaces.

Chern's conjecture concerns the same curvature condition. For each fixed
dimension, it asks whether the possible constant values of $S$ form a
discrete set \cite{Chern,CDK}. A stronger conjecture asks whether every
closed minimal hypersurface with constant $S$ is isoparametric, meaning
that all its principal curvatures are constant. Tang and Yan proved
$\lambda_1=n$ for minimal isoparametric hypersurfaces. They observed that
the stronger conjecture would therefore give the eigenvalue equality in
the constant-$S$ case \cite[Theorem~1.1 and Remark~1.3]{TangYan}.
Here we use the constancy of $S$ directly in the energy estimate.

For a general embedded minimal hypersurface, Choi and Wang extended a
first eigenfunction harmonically to one side of $M$. Reilly's identity
then gave the lower bound \cite{ChoiWang,Reilly}
\[
  \lambda_1\ge\frac{n}{2}.
\]
Duncan, Sire and
Spruck improved this bound by a positive quantity depending on $n$ and
$\max_M|A|$ \cite{DSS}. Other estimates appear in \cite{JTZ,Yu}.
Zhao used both the maximum and the minimum of $|A|$. When $S$ is constant,
his estimate gives \cite[Theorem~1.2]{Zhao}
\[
  \lambda_1>\left(\frac{1}{2}+c\right)n
\]
for an absolute constant $c>0$. For certain surfaces in $\Sph^3$, symmetry determines $\lambda_1$. Choe and Soret proved $\lambda_1=2$ for Lawson's embedded
surfaces $\xi_{m,k}$ and for the Karcher--Pinkall--Sterling surfaces
\cite[Corollary~3.3]{ChoeSoret}. Zeng has proposed a variational proof
of the eigenvalue equality for general embedded minimal hypersurfaces
\cite{Zeng}.

Several curvature results help reduce our proof to one range of $S$.
Simons' identity excludes $0<S<n$ \cite{Simons}. When $S=n$, the
hypersurface is a minimal Clifford product \cite{CDK,Lawson}. Chang
classified the closed hypersurfaces of $\Sph^4$ with constant mean
and scalar curvatures \cite{Chang}. In higher dimensions, Peng and Terng
proved that there is a gap above $S=n$ \cite{PengTerng}. Suh and Yang
proved the implication \cite{SuhYang}
\[
  S>n \quad\Longrightarrow\quad S>\frac{10n}{7}.
\]
For $n\ge4$ and $S>n$, recent preprints give stronger bounds.
Chen obtains $S>\frac{85n}{59}$ \cite{Chen}, and Cheng, Li, Wei and Yin
obtain $S>\frac{3n}{2}$ \cite{CLWY}. The bound of Suh and Yang is enough
for our proof.

\subsection{The single-layer potential and its Dirichlet energy}
Suppose that $\lambda_1<n$ or that the coordinate functions do not span
the eigenspace at $n$. In either case there is a nonzero eigenfunction $f$
with eigenvalue $0<\lambda\le n$ such that $f$ is orthogonal to every
coordinate function. We use $f$ to construct a function $u$
on the sphere. The hypersurface $M$ divides the sphere into two domains.
The function $u$ is harmonic in each domain and continuous across $M$.
We prescribe the jump in its normal derivative to be a constant multiple
of $f$. These conditions determine $u$ up to an additive constant
\cite[Theorem~A.1]{IK}. The potential has density proportional to $f$.
Its value on $M$ is determined by solving the equation. Choi and Wang,
and later Zhao, prescribed the eigenfunction itself as the boundary value
of their harmonic extensions \cite{ChoiWang,Zhao}.

The proof compares two estimates for the Dirichlet energy of $u$.
For the lower bound, we apply Reilly's identity on both sides of $M$.
The two restrictions of $u$ agree on $M$, and the two outward unit
normals are opposite. Thus the terms involving the second fundamental
form cancel when we add the two identities. The remaining boundary
integral contains $f$. Integration by parts on $M$ allows us to use
$Lf=\lambda f$ in that integral. Section~\ref{sec:reduction} gives the
resulting identity for the Hessian of $u$.

Section~\ref{sec:flux} then extends $u$ to Euclidean space by homogeneity
of degree one, separately on the cones over the two domains. Along $M$,
the difference of the two Euclidean gradients is proportional to $f\nu$,
where $\nu$ is a unit normal to $M$. Differentiating $f\nu$ in a
direction tangent to $M$ gives one term involving $df$ and another
involving the shape operator $A$. We integrate along the normal geodesics
on each side of $M$ to obtain the lower bound. Here the hypothesis on $S$ gives
\[
  \diver(A^2)=\frac{\nabla S}{2}=0.
\]

For the upper bound, we expand $u$ in spherical harmonics. The terms of
degrees zero and one vanish. We estimate the terms of degrees two and
three together using the identity
\[
 (n+1)\mathbb E(XX^{\mathsf T})+\mathbb E(\nu\nu^{\mathsf T})=I,
\]
where $\mathbb E$ denotes the average over $M$. For the higher degrees
we use the heat kernel on the sphere at time $\frac{1}{6S}$. After these
estimates it is enough to compare the bounds at $\lambda=n$
(Sections~\ref{sec:heats}--\ref{sec:green-corners}). At this value
the lower bound exceeds the upper bound, so no such eigenfunction exists.

Section~\ref{sec:principal-spectrum} uses the principal curvatures to bound the
eigenvalues of functions orthogonal to $V_X$. It also gives a condition for
$\lambda_1=n$ when $S$ varies.

\section{Reduction to one eigenfunction and its potential}\label{sec:reduction}
Write $\bar g$ and $\bar\nabla$ for the round metric and its connection.
We use $\nabla$ for the connection on $M$. Let $\nu$ be a unit normal
to $M$. We use the Hilbert--Schmidt norm for $|A|$ and the sign convention
\[
 A(v)=-\bar\nabla_v\nu.
\]

We first compute the dimension of the space $V_X$ spanned by the
coordinate functions. Suppose that $M$ is an equator. Choose orthonormal
coordinates in $\mathbb R^{n+2}$ so that $X^{n+2}=0$ on $M$. The other
$n+1$ coordinate functions are linearly independent on this equator.
Hence $\dim V_X=n+1$.

Now suppose that $M$ is not an equator. If the coordinate functions were
linearly dependent on $M$, there would be a nonzero $a\in\mathbb R^{n+2}$
such that $\langle X,a\rangle=0$ on $M$. Then $M$ would be contained in
the great sphere $a^\perp\cap\Sph^{n+1}$. Both $M$ and this great sphere
have dimension $n$, so the inclusion is a local diffeomorphism. Its image is open
in the great sphere and is also closed because $M$ is compact. Since the
great sphere is connected, $M$ would equal the great sphere. This
contradicts our assumption. Thus the $n+2$ coordinate functions are
linearly independent on $M$, and $\dim V_X=n+2$.

\paragraph{The known cases.}
On an equator the theorem follows from the spectrum of the unit
$n$-sphere. This also covers $n=1$, since every closed embedded
geodesic in $\Sph^2$ is a great circle.
If $S=n$, the rigidity theorems of Chern, do Carmo and Kobayashi
and of Lawson \cite{CDK,Lawson} identify $M$ with the product
\[
 \Sph^k\!\left(\sqrt{\frac{k}{n}}\right)
 \times \Sph^{n-k}\!\left(\sqrt{\frac{n-k}{n}}\right).
\]
Each factor has first positive eigenvalue $n$. Its coordinate functions
span that eigenspace, and its next distinct eigenvalue is greater than $2n$.
On the product, the eigenvalues are sums of eigenvalues of the two
factors, and products of eigenfunctions span the corresponding
eigenspaces. Only $0+n$ and $n+0$ give the eigenvalue $n$. Thus the
coordinate functions from the two factors span its eigenspace. The
next distinct eigenvalue is $2n$.

In dimension two the principal curvatures are $\pm\sqrt{\frac{S}{2}}$.
When $S$ is constant, the classification gives the equator and the
minimal Clifford torus \cite{Cecil}. In dimension three the only
additional case is the minimal Cartan hypersurface \cite{Chang,Cecil}.
The first two distinct positive eigenvalues of the Cartan hypersurface
are $3$ and $6$, and the eigenspace at $3$ has dimension five
\cite[Proposition~5.1]{MOU}. This eigenspace contains $V_X$. Since the
Cartan hypersurface is not an equator, $\dim V_X=5$, so the two spaces
are equal.

For constant $S$, Simons' identity \cite{Simons} gives
\begin{equation}\label{sim}
                       |\nabla A|^2=S(S-n).
\end{equation}
This rules out $0<S<n$, and $S=0$ gives an equator. We may therefore
assume $n\ge4$ and $S>n$. By the pinching theorem of Suh and Yang
\cite{SuhYang}, our two parameters satisfy
\begin{equation}\label{eq:domain}
            0<x=\frac{n}{S}<\frac{7}{10},\qquad 0<\rho=\frac{1}{n}\le\frac{1}{4}.
\end{equation}

\paragraph{The eigenfunction to exclude.}
Write $d=n+1$. Let $dx$ be volume on $\Sph^d$ and let $d\sigma$ be
area on $M$. We use the rescaled measures
\[
 d\pi=\frac{dx}{|\Sph^d|},\qquad d\mu=\frac{d\sigma}{\sqrt S\,|\Sph^d|}.
\]
Since $V_X$ contains nonzero functions and $Lh=nh$ for every $h\in V_X$,
we have $\lambda_1\le n$. Suppose that one of the two spectral conclusions
of Theorem~\ref{main} fails. If $\lambda_1<n$, self-adjointness of $L$
implies that every first eigenfunction is orthogonal to $V_X$.
If $\lambda_1=n$, then $V_X$ must be a proper subspace of $\ker(L-n)$.
In this case we choose a nonzero function in $\ker(L-n)$ orthogonal to
$V_X$. In both cases we obtain a real eigenfunction $f$ with
\[
 Lf=\lambda f,\qquad 0<\lambda\le n,\qquad f\perp V_X.
\]
Since $\lambda>0$, $f$ is also orthogonal to the constants. After rescaling $f$, these conditions read
\[
 \int_M f\,d\mu=0,\quad \int_M fX\,d\mu=0,\quad \int_M f^2\,d\mu=1 .
\]
Put $p=\frac{\lambda}{n}$.

\paragraph{The single-layer potential.}
The hypersurface $M$ separates the sphere into two domains, $\Omega_+$
and $\Omega_-$. Choose $\nu$ to point from $\Omega_-$ into $\Omega_+$.
The outward unit normal is $\nu_+=-\nu$ for $\Omega_+$ and
$\nu_-=\nu$ for $\Omega_-$.

Put $\alpha=\frac{n}{\sqrt S}$. We seek a continuous function $u$ on the
sphere. Write $u_\pm$ for the restriction of $u$ to $\Omega_\pm$ and
$g$ for the common trace on $M$. We require $u_+$ and $u_-$ to be
harmonic in their respective domains and to satisfy
\[
 u_+|_M=u_-|_M=g,\qquad
 \partial_{\nu_+}u_++\partial_{\nu_-}u_-=\alpha f .
\]
With both derivatives taken in the direction $\nu$, the second
condition reads
\[
 \partial_\nu u_- - \partial_\nu u_+ = \alpha f.
\]
The prescribed function is $f$. The boundary value $g$ is an unknown. Since
$\int_M f\,d\sigma=0$, the single-layer theorem gives a solution that
is smooth up to $M$ on each side \cite[Theorem~A.1]{IK}.
Any two solutions differ by a constant. We choose $u$ to have mean zero.
The normalizations of $u$ and $f$ give
\[
 \int_{\Sph^d}u\,d\pi=0,\qquad
 N_f=\int_M f^2\,d\sigma=\sqrt S\,|\Sph^d|.
\]

Let $E(u)$ denote the total Dirichlet energy of $u$. Green's formula on
the two domains gives
\[
 \begin{aligned}
 E(u)
 &=\sum_\pm\int_{\Omega_\pm}|\bar\nabla u_\pm|^2\,dx \\
 &=\int_M g\bigl(\partial_{\nu_+}u_++\partial_{\nu_-}u_-\bigr)\,d\sigma
 =\alpha\int_M gf\,d\sigma.
 \end{aligned}
\]
If $E(u)=0$, both $u_+$ and $u_-$ are constant and both normal
derivatives vanish on $M$. This contradicts the boundary condition
because $\alpha f$ is nonzero. Hence $E(u)>0$.

Apply Reilly's identity \cite[Theorem~1]{Reilly} to $u_+$ and $u_-$ and
add the results. The terms involving mean curvature vanish by minimality.
The terms involving the second fundamental form cancel because the two
boundary values agree and the outward normals are opposite. The Ricci
tensor of the sphere is $n\bar g$. Thus
\begin{equation}\label{reilly}
 \sum_\pm\int_{\Omega_\pm}|\bar\nabla^2 u_\pm|^2\,dx+nE(u)
       =2\alpha\int_M f Lg\,d\sigma .
\end{equation}
Since $L$ is self-adjoint and $Lf=\lambda f$, we have
\[
 \int_M fLg\,d\sigma=\lambda\int_M fg\,d\sigma.
\]
The energy formula above now gives
\begin{equation}\label{layer:energy}
 E_f:=\frac{E(u)}{n|\Sph^d|}>0,\qquad
 \sum_\pm\int_{\Omega_\pm}|\bar\nabla^2u_\pm|^2\,dx=(2\lambda-n)E(u).
\end{equation}
In particular, $2\lambda\ge n$. Equality would make the Hessian
vanish on both sides, so each gradient would be parallel. On a round
sphere a parallel vector field on an open set vanishes there because
the Ricci tensor is positive. Then both normal derivatives would vanish, contradicting $\alpha f\ne0$. We conclude that
\[
 \frac{1}{2}<p\le1.
\]
The rest of the proof gives a lower and an upper
bound for the same energy $E_f$.

\subsection{Spherical harmonics of the Green potential}
\label{sec:scalar-green}

Let $\mathcal H_j$ be the space of real spherical harmonics of
degree $j$ on $\Sph^d$. These spaces give the orthogonal decomposition
of $L^2(d\pi)$ \cite[Chap.~5]{ABR}. The Euclidean Laplacian in polar
coordinates gives the eigenvalue of $-\bar\Delta$ on $\mathcal H_j$ as
\[
 j(j+d-1)=n\sigma_j,\qquad \sigma_j=j(1+j\rho).
\]
For a finite measure $\eta$ on $\Sph^d$, let $\Pi_j\eta\in\mathcal H_j$
denote its degree-$j$ spherical harmonic component, defined by
\[
 \langle \Pi_j\eta,\varphi\rangle_{L^2(d\pi)}
 =\int_{\Sph^d}\varphi\,d\eta
 \qquad(\varphi\in\mathcal H_j).
\]
We regard $f\,d\mu$ as a measure on the sphere supported on $M$ and set
\[
 a_j=\bigl\|\Pi_j(f\,d\mu)\bigr\|_{L^2(d\pi)}^2.
\]
The harmonics of degrees zero
and one are the constants and the restrictions of linear functions,
respectively. The orthogonality conditions on $f$ thus give $a_0=a_1=0$.

Let $H_s$ denote the heat semigroup on $\Sph^d$. Write $K_s$ for its
kernel with respect to $d\pi$. For every finite measure $\eta$,
\[
 H_s\eta(a)=\int K_s(a,b)\,d\eta(b).
\]
For $w>0$, set $s=\frac{w}{n}$ and define
\[
 B(w)=\sum_{j\ge2}e^{-w\sigma_j}a_j,\qquad
 c(w)=\bigl\|H_s\mu\bigr\|_\infty.
\]
The heat kernel expansion gives
\[
 B(w)=\iint K_s(a,b)f(a)f(b)\,d\mu(a)d\mu(b).
\]
Indeed, the series is the squared $L^2(d\pi)$ norm of
$H_{\frac{s}{2}}(f\,d\mu)$, which is smooth for $s>0$.

\begin{lemma}[The Green energy in spherical harmonics]\label{lem:green-trace}
For every $t>0$,
\begin{equation}\label{green}
 \begin{aligned}
 E_f &= \sum_{j\ge2}\frac{a_j}{\sigma_j} \\
 &\le \int_0^t B(w)\,dw+\frac{B(t)}{\sigma_4}
   +\sum_{j=2,3}e^{-t\sigma_j}
             \left(\frac{1}{\sigma_j}-\frac{1}{\sigma_4}\right)a_j.
 \end{aligned}
\end{equation}
\end{lemma}
\begin{proof}
Apply Green's formula on each side and add the boundary integrals. The
prescribed sum of the normal derivatives gives
\[
 \frac1n\int_{\Sph^d}\langle du,d\varphi\rangle\,d\pi
     =\int_M f\varphi\,d\mu
\]
for every smooth $\varphi$ on the sphere. If $\varphi\in\mathcal H_j$,
integration by parts on the left gives
\[
 \sigma_j\langle u,\varphi\rangle_{L^2(d\pi)}
 =\langle \Pi_j(f\,d\mu),\varphi\rangle_{L^2(d\pi)}.
\]
Thus the degree-$j$ part of $u$ is $\frac{\Pi_j(f\,d\mu)}{\sigma_j}$ for $j\ge1$.
The two smooth restrictions of $u$ have the same trace on $M$, so together
they define a function in $H^1(\Sph^d)$. Parseval's identity for the
Dirichlet form now gives the equality in \eqref{green}.

By Tonelli's theorem, termwise integration of the nonnegative series for
$B(w)$ gives
\[
 E_f=\int_0^\infty B(w)\,dw,\qquad
 \int_t^\infty B(w)\,dw
 =\sum_{j\ge2}e^{-t\sigma_j}\frac{a_j}{\sigma_j}.
\]
In the sum on the right we keep the terms with $j=2$ and $j=3$ and use
$\frac{1}{\sigma_j}\le\frac{1}{\sigma_4}$ for each $j\ge4$. The resulting upper bound
is \eqref{green}.
\end{proof}

\section{A flux lower bound for the Green energy}\label{sec:flux}

We first estimate the energy of $u$ from below. For $y$ on either side of
$M$ and $r>0$, set $\widetilde u(ry)=ru(y)$. Its Euclidean gradient at $y$ is
\[
 F(y)=\nabla^{\mathbb R}\widetilde u(y)=\bar\nabla u(y)+u(y)y .
\]
The Euclidean gradient of $\widetilde u$ is constant on each ray from the origin.
Thus the radial vector lies in the kernel of its Euclidean Hessian.
On the tangent space to the sphere, the Hessian is the symmetric tensor
$dF=\bar\nabla^2u+u\bar g$.
Let $F_-$ and $F_+$ denote the traces of $F$ from $\Omega_-$ and $\Omega_+$, respectively.
The traces $u_-$ and $u_+$ agree on $M$, so their tangential derivatives also agree.
The prescribed jump in the normal derivative gives $F_--F_+=\alpha f\nu$.
We will differentiate this identity along $M$.

We also need an upper bound for $\int|dF|^2\,dx$.
The cross term in $|\bar\nabla^2 u+u\bar g|^2$ is
$2u\bar\Delta u=0$ on both sides. By Section~\ref{sec:scalar-green},
$u$ is orthogonal to the spherical harmonics of degrees zero and one.
The eigenvalue in degree two is $2(n+2)$. Hence
\[
 \int u^2\,dx\le \frac{E(u)}{2(n+2)}.
\]
Equation~\eqref{layer:energy} now gives
\begin{equation}\label{flat:energy}
 \int |dF|^2\,dx
 \le n(2p-1+\varepsilon)E(u),\qquad
 \varepsilon=\frac{n+1}{2n(n+2)}.
\end{equation}
The integral on the left is the sum of the integrals over the two open domains.

\subsection{The normal variation}
Use the normal $\nu$ pointing into $\Omega_+$ and put
\begin{align*}
 X_\theta&=\cos\theta X+\sin\theta\nu, & \nu_\theta&=\partial_\theta X_\theta,\\
 P_{\theta}&=\cos\theta I-\sin\theta A, & J_{\theta}&=\det P_{\theta}.
\end{align*}
The map $X_\theta$ moves each point of $M$ through signed distance
$\theta$ along its normal geodesic. Its differential on $TM$ is $P_{\theta}$,
whose eigenvalues are $\cos\theta-\kappa_i\sin\theta$.
We next determine a collar on each side of $M$. The closures of
$\Omega_-$ and $\Omega_+$ are complete. Both have smooth compact minimal
boundary $M$ and Ricci tensor $n\bar g$. Each exceptional case in
Howard's theorem has a direction of zero Ricci curvature. Thus the theorem
gives equality of the minimum inward cut distance and the minimum focal
distance in each domain \cite[Theorem~3]{Howard}. A focal point can occur only
when an eigenvalue of $P_{\theta}$ vanishes. Minimality gives $\sum_i\kappa_i=0$,
so Cauchy--Schwarz gives
\[
 \max_i\kappa_i^2\le\left(1-\frac{1}{n}\right)S.
\]
At $|\theta|=\arctan(\frac{1}{\sqrt S})$ we still have
$|\kappa_i\tan\theta|\le\sqrt{1-\frac1n}<1$. The focal radius on
each side is therefore strictly greater than this distance.
It follows that the normal coordinates are one-to-one for
$|\theta|\le\arctan\!\left(\frac{1}{\sqrt S}\right)$. The positive and negative values of
$\theta$ give collars on opposite sides of $M$.

Write $b=\frac{A}{\sqrt S}$ for the normalized shape operator. On each
normal ray we use the variables
\[
 z=\sqrt S\,|\tan\theta|,\qquad y=1-z,\qquad
 s=\operatorname{sgn}\theta,\qquad c=\cos\theta.
\]
Choose three radial polynomials $H_1,H_2,H_3$ with $H_1(0)=1$
and $H_2(0)=0$. Define
\[
 h_j=y^{\frac{3}{2}}H_j\quad (j=1,2,3),\qquad D=\sqrt S(h_1b+sh_2b^2).
\]
Thus $h_1$ and $sh_2$ are the coefficients of $b$ and $b^2$ in
$\frac{D}{\sqrt S}$. We use $h_3$ for the normal component of the test form.
The factor $y^{\frac{3}{2}}$ makes each squared radial derivative divisible
by $y$. We will use this factor to bound the terms containing $\frac{1}{J_{\theta}}$
near the outer edge of the collar. Define a vector-valued one-form
on $M$ by
\[
         \Phi_\theta(v)=fDv+h_3\,df(v)\nu .
\]
Set $F_\theta=F\circ X_\theta$.
We take inner products of one-forms using the metric on $M$ at $\theta=0$.
The dot denotes the Euclidean inner product. For two vector-valued
one-forms, the colon denotes the sum of these inner products over an
orthonormal basis of $TM$.
The product rule and integration by parts on the closed manifold $M$ give
\begin{equation}\label{flux}
 \frac{d}{d\theta}\int_M d_MF_\theta:\Phi_\theta\,d\sigma
 =\int_M\left(d_MF_\theta:\partial_\theta\Phi_\theta
   -\partial_\theta F_\theta\cdot
                     \operatorname{div}^{\mathbb R}_M\Phi_\theta\right)d\sigma .
\end{equation}
Here we use the Euclidean connection to differentiate the values of
$\Phi_\theta$. Since $h_2(0)=0$, $\Phi_\theta$ has a common trace at
$\theta=0$. Integrate \eqref{flux} separately on each side,
in the direction of increasing $\theta$. At the two outer ends,
$y=0$, so $\Phi_\theta=0$. The boundary contribution at $\theta=0$ is the integral of
$(d_MF_--d_MF_+):\Phi_0$ over $M$. Since $d_M(\alpha f\nu)=\alpha(df\otimes\nu-fA)$,
this boundary term is
\[
       -\alpha\bigl(S-\lambda h_3(0)\bigr)N_f .
\]

We now write the integral on the right of \eqref{flux} using the
volume measure $dx$ on the sphere. Write $M_\theta=X_\theta(M)$.
The tangent spaces to $M$ and $M_\theta$ at corresponding points
are the same Euclidean subspace, perpendicular to both $X$ and $\nu$.
Define the vector-valued one-form $T$ on the collar by
\[
 T(v)=J_{\theta}^{-1}(\partial_\theta\Phi_\theta)(P_{\theta}v)\quad (v\in TM_\theta),
 \qquad
 T(\nu_\theta)=-J_{\theta}^{-1}\operatorname{div}^{\mathbb R}_M\Phi_\theta .
\]
The chain rule gives $d_MF_\theta=dF\circ P_{\theta}$ and
$\partial_\theta F_\theta=dF(\nu_\theta)$. The volume
element is $dx=J_{\theta}\,d\sigma\,d\theta$. Thus the right-hand side of
\eqref{flux}, integrated over both collars, is $\int dF:T\,dx$.
Project the values of $T$ onto $T_{X_\theta}\Sph^{n+1}$.
Let $T_{\rm s}$ denote the symmetric part of the resulting endomorphism. Since $dF$ is tangent to the sphere and
symmetric, $dF:T=dF:T_{\rm s}$.
Cauchy--Schwarz and \eqref{flat:energy} give
\begin{equation}\label{combined}
 \alpha^2(S-\lambda h_3(0))^2 N_f^2
 \le n(2p-1+\varepsilon)E(u)I_D,\qquad
                 I_D=\int|T_{\rm s}|^2\,dx .
\end{equation}

To estimate $I_D$, split $T_{X_\theta}\Sph^{n+1}$ into
$TM_\theta$ and the line spanned by $\nu_\theta$. We compute the
tangential, mixed, and normal parts of $T_{\rm s}$.
Let $b_i=\frac{\kappa_i}{\sqrt S}$ be the eigenvalues of the normalized
shape operator $b$. For $k\ge1$, define the normalized power sum
\[
 f_k=\sum_{i=1}^n b_i^k=\tr(b^k)=\frac{\tr(A^k)}{(\sqrt S)^k}.
\]
Minimality gives $f_1=0$, and the definition of $S$ gives $f_2=1$. The Codazzi equation and minimality give
$\operatorname{div}_M A=0$ and $\operatorname{div}_M A^2=\frac{1}{2}\nabla S$.
Since $S$ is constant, $\operatorname{div}_M A^2=0$. The Euclidean Gauss
formula then gives
\[
 \operatorname{div}^{\mathbb R}_M\Phi_\theta=
 (D-h_3 A)\nabla f+
 f\bigl((Sh_1+sSh_2f_3-\lambda h_3)\nu-s\sqrt S\,h_2 X\bigr).
\]
At $X_\theta$ the projections of $\nu$ and $X$ onto the tangent
space of the sphere are $c\nu_\theta$ and $-\sin\theta\,\nu_\theta$, respectively.
The tangential block of $T_{\rm s}$ is therefore $\frac{fP_{\theta}\,\partial_\theta D}{J_{\theta}}$.
The mixed block is $\operatorname{Sym}(\nu_\theta\otimes V)$, where
$\operatorname{Sym}$ denotes the symmetric part. The vector $V$ and
the normal entry are
\begin{align*}
 V&=\frac{\bigl(c(\partial_\theta h_3)P_{\theta}-(D-h_3 A)\bigr)\nabla f}{J_{\theta}},\\
 T_{\rm s}(\nu_\theta,\nu_\theta)
 &=-\frac{cSf}{J_{\theta}}\left(h_1+sh_2f_3+\frac{zh_2}{S}-pxh_3\right),
 \qquad x=\frac{n}{S}.
\end{align*}

\subsection{A choice of the collar fields}
With $x=\frac{n}{S}$ and $y=1-z$, choose the three profiles
\begin{equation}\label{new:profile}
 H_1=1-\frac{zy}{3},\qquad H_2=\frac{1-y^2}{2},\qquad H_3=2H_2-\frac{3}{2}.
\end{equation}
We record three properties of these polynomials.
First, $H_2'=y$ and $H_2(0)=0$, where the prime denotes the derivative
with respect to $z$. Second, $H_3=-\frac{1}{2}-y^2$, so $h_3(0)=-\frac{3}{2}$.
These identities give
\[
 \frac{S-\lambda h_3(0)}{S}=1+\frac{3px}{2},
 \qquad
 yH_3'-\frac{3H_3}{2}=\frac{3}{4}+\frac{7y^2}{2}>0,
\]
where the first identity is used in \eqref{combined} and the inequality
holds throughout the collar.
Third, $zH_1-H_2=z^2\left(\frac{1}{2}-\frac{y}{3}\right)\ge0$. Hence $0\le H_2\le zH_1\le H_1$.
The bound $H_2\le zH_1$ will control the sum of the integrals on the two normal rays.
\begin{proposition}\label{prop:naturalflux}
Put $x_*=\frac{7}{10}$, $\tau(x)=1+\frac{3x}{2}$ and let $\ell$ be the affine function with
\[
              \ell(0)=\frac{4}{13},\qquad \ell(x_*)=\frac{6\tau(x_*)}{25}.
\]
For the normalized eigenfunction $f\perp V_X$ of
Section~\ref{sec:reduction}, assume \eqref{eq:domain} and $\frac{1}{2}<p\le1$.
Then
\begin{equation}\label{eq:naturalflux}
 E_f^{-1}\le \frac{2p-1+\varepsilon}{\ell(x)}
          \left(\frac{2\tau(x)}{2+3px}\right)^2 .
\end{equation}
\end{proposition}

We will prove the estimate for $0\le x\le\frac{7}{10}$ and $0\le\rho\le\frac{1}{4}$.
We first bound the integrals over the two rays for fixed principal
curvatures. Convexity then reduces the estimate to $x=0$ and $x=\frac{7}{10}$.

\subsection{The paired normal Jacobians}\label{sec:collar-normal}
Fix a point of $M$ and consider its normalized principal curvatures
$b_1,\ldots,b_n$. Up to this point we have used one fixed normal field.
We can now choose how to label the two rays at each point.
Reversing the signs of $\nu,\theta,A$ leaves $X_\theta,P_{\theta},J_{\theta}$ unchanged
and changes the signs of both $\nu_\theta$ and $\Phi_\theta$.
The squared norm $|T_{\rm s}|^2$ in \eqref{combined} is therefore
independent of the choice of labels. At each point we label the two rays
so that $f_3\ge0$. This choice is made in the scalar integrand after
all derivatives have been taken.

Put $t=f_3^2$. We need three inequalities for the power sums $f_3,f_4,f_5,f_6$.
Since $\sum_i b_i^2=1$, we may give each index $i$ probability $b_i^2$.
The random variable with value $b_i$ at index $i$ has mean $f_3$ and
second moment $f_4$. We claim that
\begin{equation}\label{cone}
    t\le f_4\le\frac{1+t}{2},\qquad 2f_6\le f_4+t,\qquad
    f_3f_5\le t+\tfrac14\left(\frac{1+t}{2}-f_4\right).
\end{equation}
The variance is $f_4-t=\sum_{i<j}b_i^2b_j^2(b_i-b_j)^2$.
Since $(b_i-b_j)^2\le2$, we get $0\le f_4-t\le1-f_4$.
For the second bound, expand
\[
       f_4+t-2f_6=\sum_{i<j}b_i^2b_j^2(b_i+b_j)^2 .
\]
The variance of $b_i^2-f_3b_i$ is nonnegative. Expanding it gives
\[
 2f_3(f_5-f_3f_4)\le t(f_4-t)+f_6-f_4^2.
\]
Inserting the first two bounds yields
\[
 2(f_3f_5-t)\le(f_4-t)\left(\frac{1}{2}-f_4\right).
\]
Here we used $t(2f_4-1-t)\le0$, which follows from the first bound in
\eqref{cone}. For $f_4\ge\frac{1}{2}$, the right-hand side is nonpositive.
For $f_4\le\frac{1}{2}$ we have $f_4-t\le\frac{1}{2}$ and $\frac{1}{2}-f_4\le\frac{1+t}{2}-f_4$,
which proves \eqref{cone}.

Let $w_x(z)=\exp\!\left(\frac{xz^2}{2}\right)$ and introduce the two radial weights
\[
 B_+=e^{-z},\qquad B_-=\frac{1-z}{1+z}e^z,\qquad
 \omega_s=\tfrac12(1+sf_3)B_++\tfrac12(1-sf_3)B_- .
\]
We claim that $\frac{y}{J_{\theta}}\le w_x\omega_s$. The factor $c^{-n}$ in $J_{\theta}^{-1}$ is
$\left(1+\frac{z^2}{S}\right)^{\frac{n}{2}}\le w_x$. To bound the remaining factors, use the identity
\[
 -\log(1-az)=az+a^2\int_0^1\frac{z^2v}{1-avz}\,dv .
\]
For fixed $z<1$, the integral is convex in $a\in[-1,1]$.
It therefore lies below the line joining its values at $a=-1$ and $a=1$.
Evaluate this bound at $a=sb_i$, multiply by $b_i^2$, and sum over $i$.
The linear terms $sb_i z$ sum to zero. Exponentiating gives
\[
 y\prod_i(1-szb_i)^{-1}
 \le \sqrt{\frac{1-z}{1+z}}\exp\!\bigl(sf_3(\operatorname{artanh}z-z)\bigr).
\]
For fixed $z$, the right-hand side is convex as a function of $sf_3$.
Its value at $sf_3=-1$ is $B_-$, and its value at $sf_3=1$ is $B_+$.
The line through these two values gives the bound $\frac{y}{J_{\theta}}\le w_x\omega_s$.
We also have $\omega_+\ge \omega_-$ because $\operatorname{artanh}z\ge z$ and
$f_3\ge0$. For a radial function $R$, use the following integrals against the
sum and difference of the weights, respectively.
\[
   L_+(R)=\int_0^1(B_++B_-)R\,dz,\qquad
   L_-(R)=\int_0^1(B_+-B_-)R\,dz .
\]

\subsection{The three radial blocks}\label{sec:radial-bulk}
Take a principal unit vector $e_i$. The eigenvalue of $P_{\theta}$ in this
direction is $c(1-szb_i)$ and $\partial_\theta z=\frac{s\sqrt S}{c^2}$.
On each ray the volume element is
\[
 dx=\frac{J_{\theta}c^2}{\sqrt S}\,d\sigma\,dz,\qquad 0\le z\le1.
\] To compute $\frac{\sqrt S\,I_D}{S^2}$, we may therefore multiply
each block of $T_{\rm s}$ by $\frac{cJ_{\theta}}{S}$ and integrate its squared
norm against $\frac{1}{J_{\theta}}\,d\sigma\,dz$. For each of the three profiles, set
\[
 d_j=yH_j'-\frac{3H_j}{2}\qquad(j=1,2,3).
\]
Then $\frac{dh_j}{dz}=\sqrt y\,d_j$.

\paragraph{The tangential and normal blocks.}
The tangential and normal contributions to $\frac{\sqrt S\,I_D}{S^2}$ have the
following densities with respect to $d\sigma\,dz$, in that order.
\begin{align*}
 &\frac{yf^2}{J_{\theta}}\sum_i(1-szb_i)^2(d_1b_i+sd_2b_i^2)^2,\\
 &\frac{yf^2}{J_{\theta}}\,y^2c^4\left(H_1+sf_3H_2+\frac{zH_2}{S}-\frac{\lambda}{S}H_3\right)^2.
\end{align*}
Both densities are nonnegative, so we can replace $\frac{y}{J_{\theta}}$ by its upper
bound $w_x\omega_s$ on each ray. To estimate the normal block we add the
values on the two rays before bounding the factor $c^4$. Put $a=H_1-pxH_3$ and $k=f_3H_2$. After multiplication by $c^2$, the two
factors inside the squares become $\widetilde a\pm\widetilde k$, where
\[
 \widetilde a=c^2\left(a+\frac{zH_2}{S}\right),\qquad \widetilde k=c^2k .
\]
Since $0\le H_2\le zH_1\le za$ and $c^2=\left(1+\frac{z^2}{S}\right)^{-1}$,
we have $0\le\widetilde a\le a$ and $0\le\widetilde k\le k$.
For $A,B\ge0$ the sum
\[
 \omega_+(A+B)^2+\omega_-(A-B)^2
    =(\omega_++\omega_-)(A^2+B^2)+2(\omega_+-\omega_-)AB
\]
is nondecreasing in $A$ and $B$. We may therefore replace
$\widetilde a$ by $H_x=H_1-xH_3\ge a$ and $\widetilde k$ by $f_3H_2$ in
the sum of the two normal squares. For fixed $b_i$, we now expand
the squares in both blocks, add the two rays, and integrate in $z$.
For these fixed curvatures, the coefficient of $f^2$ is bounded by
\begin{equation}\label{P0}
 \begin{aligned}
 Q_0(x)={}&L_+(C_0)+f_4 L_+(C_4)+f_6 L_+(C_6)\\
          &+t\bigl(L_+(C_2)+L_-(C_1)\bigr)+f_3f_5\,L_-(C_5).
 \end{aligned}
\end{equation}
Each coefficient depends on $z$ and $x$. Even indices correspond to the sum
of the weights and odd indices to their difference.

\noindent The six coefficients in \eqref{P0} are
\[
\begin{aligned}
 C_0&=w_x(d_1^2+y^2H_x^2),\\
 C_1&=w_x(-2zd_1^2+2d_1d_2+2y^2 H_x H_2),\\
 C_2&=w_x\,y^2H_2^2,\\
 C_4&=w_x(z^2d_1^2+d_2^2-4zd_1d_2),\\
 C_5&=w_x(-2zd_2^2+2z^2d_1d_2),\\
 C_6&=w_x\,z^2d_2^2.
\end{aligned}
\]
For a given $x$, let $\Lambda_0$ be the supremum of $Q_0(x)$ over all tuples
$(b_1,\ldots,b_n)$ with $f_1=0$, $f_2=1$, and $f_3\ge0$. Thus $\Lambda_0$
bounds the coefficient of $f^2$ independently of the principal curvatures.

\paragraph{The mixed block.}
After multiplication by $\frac{cJ_{\theta}}{S}$ the mixed block is
$\operatorname{Sym}(\nu_\theta\otimes\widetilde V)$,
where $\widetilde V=\frac{cJ_{\theta}}{S}V$. In a principal frame,
\[
 \widetilde V_i=\frac{c\sqrt{y}\,e_i(f)}{\sqrt S}(M_0+sM_1),
 \qquad \beta=b_i,
\]
where the factors even and odd in $s$ are determined by
\[
 M_0=\bigl(y(H_3-H_1)-zd_3\bigr)\beta,
 \qquad M_1=d_3-yH_2\beta^2.
\]
Since $c\le1$ and
$|\operatorname{Sym}(\nu_\theta\otimes\widetilde V)|^2=\frac{|\widetilde V|^2}{2}$, the coefficient of $\frac{|df|^2}{S}$ after integration in $z$
is at most
\[
 \Lambda_1=\tfrac12\sup_{\beta,f_3}
    \bigl(L_+(w_x(M_0^2+M_1^2))+f_3\,L_-(2w_x\,M_0M_1)\bigr).
\]
Here $\beta$ is one of the $b_i$, and $f_3=\sum_j b_j^3$ is computed from
the same tuple. For each such pair,
\begin{equation}\label{mix:spectrum}
 0\le 1-f_3\beta\le 3(1-\beta^2).
\end{equation}
Indeed $|f_3-\beta^3|\le\sum_{j\ne i}|b_j|^3\le1-\beta^2$,
where $b_i=\beta$. Thus
$1-f_3\beta\le(1-\beta^2)(1+\beta^2+|\beta|)$.

Integrating the three blocks and using
$\int_M|df|^2\,d\sigma=\lambda N_f\le nN_f$ gives
\[
 \frac{\sqrt S\,I_D}{S^2N_f}\le \Lambda_0+x \Lambda_1.
\]
Now $E(u)=\alpha N_fE_f$, $\alpha=\frac{n}{\sqrt S}$, and $h_3(0)=-\frac{3}{2}$.
Substitution in \eqref{combined} yields
\begin{equation}\label{upbound}
 \left(1+\frac{3px}{2}\right)^2E_f^{-1}\le(2p-1+\varepsilon)(\Lambda_0+x\Lambda_1).
\end{equation}

\subsection{Proof of the flux bound}\label{sec:radial-proof}
\begin{proof}[Proof of Proposition~\ref{prop:naturalflux}]
By \eqref{upbound}, it is enough to prove
$(\Lambda_0+x\Lambda_1)\ell(x)\le\tau(x)^2$. We first bound $\Lambda_1$ by a function
of $x$ alone. We then use convexity to reduce the estimate for
$\Lambda_0+x\Lambda_1$ to the two values $x=0$ and $x=x_*=\frac{7}{10}$.

\smallskip\noindent\emph{A bound for the mixed block.}
The mixed block depends on three radial coefficients. Set
\[
 m_0=y(H_3-H_1)-zd_3,\quad m_1=d_3=\frac{3}{4}+\frac{7y^2}{2},\quad m_2=-yH_2 .
\]
Thus $M_0=m_0\beta$ and $M_1=m_1+m_2\beta^2$. The expression inside
the supremum defining $2\Lambda_1$ is the polynomial
\[
 q_0+q_2\beta^2+q_4\beta^4+f_3(q_1\beta+q_3\beta^3).
\]
The coefficient $q_j$ corresponds to the power $\beta^j$.
For $q_0,q_2,q_4$, apply $L_+$ to $w_x\,m_1^2$, $w_x(m_0^2+2m_1m_2)$
and $w_x\,m_2^2$, respectively. For $q_1,q_3$, apply $L_-$ to
$2w_x\,m_0m_1$ and $2w_x\,m_0m_2$. The even coefficients are nonnegative.
For $q_2$, use
$m_0^2\ge4yz(H_1-H_3)m_1\ge2ym_1H_2$ and $m_2=-yH_2$, where $H_2=\frac{z(1+y)}{2}$.
Also $q_1\le0\le q_3$, since $B_+\ge B_-$ and $m_0,m_2\le0\le m_1$.

We compare the polynomial at $(f_3,\beta)$ with its value at $(1,1)$.
The terms $q_4\beta^4$ and $f_3q_3\beta^3$ are at most $q_4$ and
$q_3$, respectively.
By \eqref{mix:spectrum}, the increase from the other terms is at
least $(q_2+3q_1)(1-\beta^2)$. To prove that $q_2+3q_1>0$, we use $B_+,B_-\ge y$, $w_x\le\frac{3}{2}$, and
\[
 B_+-B_-=\frac{2(z\cosh z-\sinh z)}{1+z}\le z^3 .
\]
To bound $B_+-B_-$, integrate $a\sinh a\le a^2\sinh 1$ from $0$ to
$z$ and use $2\sinh 1<3$. We obtain
\[
 q_2+3q_1\ge
  \int_0^1\left(2y(m_0^2+2m_1m_2)+9z^3m_0m_1\right)\,dz
  >0 .
\]
The strict inequality follows by integrating each monomial in $z$. Evaluating at $(1,1)$ now gives $\Lambda_1\le Q_1(x)$, where
\[
  Q_1(x)=\tfrac12\sum_{s=\pm1}\int_0^1 w_x B_s
                 \bigl(m_0+s(m_1+m_2)\bigr)^2\,dz .
\]

\smallskip\noindent\emph{Reduction to two values of $x$.}
Fix the normalized principal curvatures $b_i$. Combine their tangential
and normal contribution with the mixed estimate by setting
$Q(x)=Q_0(x)+xQ_1(x)$, with $Q_0$ as in \eqref{P0}.
All Taylor coefficients of $Q$ at
$x=0$ are nonnegative. To see this for $Q_0$, consider the squares on
the two rays before expanding them. Apart from the exponential factor
$w_x$, the tangential square is independent of $x$. The normal square
contains the factor
\[
                 \bigl(H_1+sf_3H_2-xH_3\bigr)^2 .
\]
Here $H_1+sf_3H_2\ge H_1-H_2\ge0$ and $-H_3>0$. For $Q_1$, the only factor that
depends on $x$ is $w_x$.

We need a lower bound for the coefficient of $x^2$ in $Q$.
Since $B_+,B_-\ge y$, the normal square contributes at least
$2\int_0^1y^3H_3^2\,dz=\frac{17}{24}$. The contribution from $xQ_1(x)$ is $Q_1'(0)$.
Keeping the term with weight $B_-$ gives
\begin{align*}
 Q_1'(0)&\ge\frac14\int_0^1 z^2y(m_0-m_1-m_2)^2\,dz\\
 &\ge 3\left(\int_0^1 z^2y(m_1+m_2-m_0)\,dz\right)^2
 =3\left(\frac{8}{35}\right)^2 .
\end{align*}
The second inequality uses Cauchy--Schwarz and $\int_0^1z^2y\,dz=\frac{1}{12}$.
The last equality follows by substituting $m_0,m_1,m_2$.
Since $\frac{17}{24}+3\left(\frac{8}{35}\right)^2>\frac{6}{7}$, the function $Q(x)-\frac{6}{7}x^2$ is convex.

We will prove the endpoint bounds
\begin{equation}\label{natural:anchors}
  Q(0)<A_0:=\tfrac{13}{4},\qquad
  Q(x_*)<A_*:=\tfrac{25}{6}\tau(x_*).
\end{equation}
Assume \eqref{natural:anchors} for the moment. Set $a=\frac{x}{x_*}$ and
$\overline A=(1-a)A_0+a A_*$. The convex function $Q(x)-\frac{6}{7}x^2$
lies below the line joining its endpoint values. Thus
\[
 Q(x)\le\overline A-\frac{6}{7}x_*^2a(1-a).
\]
From the definition of
$\ell$, we have
\[
 \overline A\,\ell(x)-\tau(x)^2
   =a(1-a)\,\ell(0)\,\frac{\bigl(A_*-A_0\tau(x_*)\bigr)^2}{A_*}.
\]
Substituting the constants in \eqref{natural:anchors} gives
$\bigl(A_*-A_0\tau(x_*)\bigr)^2<\frac{6}{7}x_*^2A_*$. Since $\ell(x)\ge\ell(0)>0$,
these inequalities give $Q(x)\ell(x)\le\tau(x)^2$ for $0\le x\le x_*$.

\smallskip\noindent\emph{Proof of the endpoint bounds.}
We use the same polynomial estimate for both weights. Let
\[
 P_4(y)=1+y+\frac{y^2}{2}+\frac{y^3}{6}+\frac{y^4}{16}.
\]
Then $e^y\le P_4(y)$ for $y\le\frac{3}{2}$.
For $y\le0$, use Taylor's theorem at zero. For $0\le y\le\frac{3}{2}$,
the tail of the exponential series starting at $\frac{y^4}{24}$ is at most
\[
 \frac{y^4}{24\left(1-\frac{y}{5}\right)}\le \frac{y^4}{16}.
\]
Set
\[
 y_\pm=\frac{xz^2}{2}\mp z,\qquad v=\frac{zy}{2},\qquad
 R_+=P_4(y_+),\quad R_-=yP_4(y_-)\left(1-\frac{z}{2}\right)(1-v+v^2).
\]
Then $w_xB_+\le R_+$ and $w_xB_-\le R_-$. The second inequality also
uses $(1+z)\left(1-\frac{z}{2}\right)=1+v$ and $\frac{1}{1+v}\le1-v+v^2$.
Apply these bounds to the squares on the two rays before expanding.
For the polynomial upper bounds, the weighted integrals of the
coefficient functions $C_j$ become
\[
 c_{j}^{+}=\int_0^1 (R_++R_-)\frac{C_j}{w_x}\,dz,\quad
 c_{j}^{-}=\int_0^1 (R_+-R_-)\frac{C_j}{w_x}\,dz .
\]
Let $\widehat Q_1$ be the integral defining $Q_1$ with $w_xB_s$ replaced
by $R_s$. Since $c_{6}^{+}\ge0$ and $2f_6\le f_4+t$, \eqref{P0} gives
\[
\begin{gathered}
 Q(x)\le D_0+D_1f_4+D_2t+D_3f_3f_5,\\[4pt]
 D_0=x\widehat Q_1+c_{0}^{+},\quad
 D_1=c_{4}^{+}+\frac{c_{6}^{+}}{2},\quad
 D_2=c_{2}^{+}+c_{1}^{-}+\frac{c_{6}^{+}}{2},\quad
 D_3=c_{5}^{-}.
\end{gathered}
\]
Recall that $t=f_3^2$. To bound the terms in $f_4$ and $f_3f_5$, set
$h=\frac{1+t}{2}-f_4$. By \eqref{cone}, $0\le h\le\frac{1-t}{2}$ and
$f_3f_5\le t+\frac{h}{4}$. Thus, if $D_3\ge0$, then
\[
 D_1f_4+D_3f_3f_5\le \frac{D_1(1+t)}{2}+D_3t+\left(\frac{D_3}{4}-D_1\right)h .
\]
Maximize first over $0\le h\le\frac{1-t}{2}$. The resulting bound is
affine in $t$, so its maximum on $[0,1]$ occurs at $0$ or $1$. This gives
\[
 Q(x)\le \max\left(D_0+\max\!\left(\frac{D_1}{2},\frac{D_3}{8}\right),\;
                      D_0+D_1+D_2+D_3\right).
\]
At each of $x=0$ and $x=\frac{7}{10}$, all the integrals defining $D_0,\ldots,D_3$
are integrals of polynomials. Substitution and the rule
\[
 \int_0^1\sum_j c_j z^j\,dz=\sum_j \frac{c_j}{j+1}
\]
give $D_3>0$ in both cases.
The displayed maximum is less than $A_0$ at $x=0$ and less than $A_*$
at $x=\frac{7}{10}$. This proves \eqref{natural:anchors}.

We have proved $Q(x)\ell(x)\le\tau(x)^2$ for every normalized set of
principal curvatures. Taking the supremum over these curvatures and
using $\Lambda_1\le Q_1(x)$ gives $(\Lambda_0+x\Lambda_1)\ell(x)\le\tau(x)^2$.
The proposition now follows from \eqref{upbound}.
\end{proof}

\section{Heat kernels on the hypersurface and the sphere}\label{sec:heats}
\subsection{Two bounds from the normal collar}
Recall that
\[
 x=\frac{n}{S},\qquad \rho=\frac{1}{n},\qquad
 c(w)=\bigl\|H_{\frac{w}{n}}\mu\bigr\|_\infty .
\]
Write $\mu_0=\mu(M)$ for the total mass and set
$D=\bigl(\tfrac12+\rho\bigr)x$. Define
\[
 I(y)=\int_0^1 e^{-yz^2}\sqrt{1-z^2}\,dz,\qquad c_I(y)=\frac1{2I(y)}.
\]
We bound $c(w)$ and the total mass $\mu_0$ by integrating along
the normal geodesics on both sides of $M$.

\begin{lemma}\label{lem:collar-cap}
Assume \eqref{eq:domain}. For every $w>0$,
\begin{equation}\label{heatcap}
 c(w)\le c_I\!\left(D+\frac{x}{4w}\right),\qquad
 \mu_0\le c_I(D).
\end{equation}
\end{lemma}
\begin{proof}
At $a\in M$, let $b_i=\tfrac{\kappa_i(a)}{\sqrt S}$ be the eigenvalues
of $b=\tfrac{A}{\sqrt S}$.
On each normal ray from $a$, let $\theta$ be the distance from $a$ and
set $z=\sqrt S\tan\theta$. For $0\le z\le1$, these coordinates cover
the two collars in Section~\ref{sec:flux}. The normalized spherical
volume element is
\[
 d\pi=\left(1+\frac{z^2}{S}\right)^{-\frac{n+2}{2}}
 \det(I\mp zb)\,d\mu\,dz .
\]
Since $\log(1+u)\le u$ for $u\ge0$ and $D=\tfrac{n+2}{2S}$, we have
\[
 \left(1+\frac{z^2}{S}\right)^{-\frac{n+2}{2}}\ge e^{-Dz^2}.
\]

Fix $\xi\in\Sph^{n+1}$ and put $s=\tfrac{w}{n}$.
Let $K_+$ and $K_-$ be the values of $K_s(\xi,\cdot)$ at distance
$\theta$ from $a$ in the directions $\nu$ and $-\nu$, respectively.
Applying Hamilton's estimate \cite[Corollary~4.4]{Hamilton} to the
smooth positive solution $K_{t+\varepsilon}(\xi,\cdot)$ at time $t=s$
and letting $\varepsilon\downarrow0$ gives
\[
 \bar\nabla^2\log K_s\ge -\frac{\bar g}{2s}.
\]
Integrate the Hessian bound twice from $a$ to distance $\theta$ on each
ray. The two initial directional derivatives of $\log K_s(\xi,\cdot)$
have opposite signs. Adding the two inequalities and exponentiating gives
\[
 \sqrt{K_+K_-}\ge K_s(\xi,a)\exp\!\left(-\frac{\theta^2}{4s}\right).
\]
The product of the two determinants is $\det(I-z^2b^2)$. We have
\[
 \sqrt{\det(I-z^2b^2)}=\sqrt{\prod_i(1-b_i^2z^2)}\ge\sqrt{1-z^2}
 \qquad(0\le z\le1),
\]
because $\sum_i b_i^2=1$ and $\prod_i(1-a_i)\ge1-\sum_i a_i$ for
$0\le a_i\le1$.

The sum of $K_+\det(I-zb)$ and $K_-\det(I+zb)$ is at least twice their geometric mean.
Insert the bounds for $\sqrt{K_+K_-}$ and $\sqrt{\det(I-z^2b^2)}$.
Use $\theta\le \tfrac{z}{\sqrt S}$ and integrate over $0\le z\le1$ and $a\in M$.
The spherical heat kernel has total mass one. Hence
\[
 1\ge\int_{\rm collars} K_s(\xi,\cdot)\,d\pi
 \ge 2 I\!\left(D+\frac{x}{4w}\right) H_s\mu(\xi).
\]
Taking the supremum over $\xi$ gives the bound for $c(w)$. Integrating the
constant function $1$ over the same collars gives $1\ge2I(D)\mu_0$.
\end{proof}

\subsection{The semicircle integral}
\begin{lemma}\label{lem:heat-cap}
On $[0,\infty)$ the function $c_I$ is increasing and concave, and $c_I^2$ is
convex. In particular,
\[
 0\le c_I'(y)\le\frac{1}{2\pi},\qquad (c_I^2)'(y)\le\frac{1}{\pi}
 \qquad (y\ge0).
\]
\end{lemma}
\begin{proof}
In the integral for $I(y)$ put $r=z^2$. After normalization the density
of $r$ on $(0,1)$ is proportional to $\sqrt{\tfrac{1-r}{r}}\,e^{-yr}$.
Write $g_y$ for this density and denote its mean and second moment by
$q=q(y)$ and $m=m(y)$, respectively. Differentiation with respect to $y$ gives
\[
 (\log c_I)'=q>0,\qquad q'=q^2-m.
\]
It follows that
\[
 \frac{c_I''}{c_I}=2q^2-m,\qquad \frac{(c_I^2)''}{2c_I^2}=3q^2-m .
\]
Thus $m\ge2q^2$ gives concavity of $c_I$, and $m\le3q^2$ gives convexity of $c_I^2$.

We compare $g_y$ with two probability densities that have the same mean
$q$. Let $g_\beta$ be the density of $r_* r_0$, where $r_0$ has density
$g_0$ and $r_*=4q$. By the beta integral, $r_0$ has mean $\tfrac14$ and
second moment $\tfrac18$. Thus $r_* r_0$ has mean $q$ and second moment $2q^2$.
Let $g_\gamma$ be the gamma density with shape $\tfrac12$ and mean $q$.
Its second moment is $3q^2$.

Suppose $y>0$. Since $q(0)=\tfrac14$ and $q'=-\operatorname{Var}(r)<0$,
we have $r_*<1$. Extend each density by zero outside its support.
On the intervals where the ratios are positive, we have
\[
\begin{aligned}
 \frac{d^2}{dr^2}\log\frac{g_y}{g_\beta}
 &=\tfrac12\bigl((r_*-r)^{-2}-(1-r)^{-2}\bigr)>0 &&(0<r<r_*),\\
 \frac{d^2}{dr^2}\log\frac{g_y}{g_\gamma}
 &=-\tfrac12(1-r)^{-2}<0 &&(0<r<1).
\end{aligned}
\]
The function $\log\!\bigl(\tfrac{g_y}{g_\beta}\bigr)$ is strictly convex and tends to
$+\infty$ as $r\uparrow r_*$. Also $g_y-g_\beta>0$ on $(r_*,1)$.
Thus $g_y-g_\beta$ changes sign at most twice and is positive after its
last sign change wherever it is nonzero. The function
$\log\!\bigl(\tfrac{g_y}{g_\gamma}\bigr)$ is strictly concave and tends to $-\infty$ as
$r\uparrow1$. Since $g_y-g_\gamma<0$ for $r>1$, the difference
$g_y-g_\gamma$ also changes sign at most twice and is negative after
its last sign change.

For either $g=g_\beta$ or $g=g_\gamma$, the two densities have equal
mass and differ on an interval. Hence $g_y-g$ changes sign at least
once. If there were a single sign change at $a$, the function
$(r-a)(g_y-g)$ would have a fixed sign and a nonzero integral.
However, equality of mass and mean gives
\[
 \int_0^\infty(r-a)(g_y-g)\,dr=0.
\]
Thus $g_y-g$ changes sign at two points $a<b$. Again using equality
of mass and mean, we obtain
\[
 \int_0^\infty r^2(g_y-g)\,dr
  =\int_0^\infty(r-a)(r-b)(g_y-g)\,dr .
\]
For $g=g_\beta$ the integrand on the right is nonnegative, so $m\ge2q^2$.
For $g=g_\gamma$ the integrand is nonpositive, so $m\le3q^2$.
By continuity both inequalities hold at $y=0$.

Since $c_I(0)=\tfrac{2}{\pi}$ and $q(0)=\tfrac14$, concavity gives
\[
 c_I'(y)\le c_I'(0)=\frac{1}{2\pi}.
\]
In $I(y)$ make the substitution $v=\sqrt y\,z$. By dominated convergence,
\[
 \lim_{y\to\infty}\sqrt{y}\,I(y)=\frac{\sqrt\pi}{2},\qquad
 \lim_{y\to\infty}\frac{c_I(y)^2}{y}=\frac1\pi .
\]
By convexity, $(c_I^2)'(y)$ is at most the slope of the chord from $y$ to any
larger argument. Sending that argument to infinity gives
$(c_I^2)'(y)\le\tfrac{1}{\pi}$.
\end{proof}

\subsection{Comparison with the heat kernel of the hypersurface}
\begin{lemma}\label{lem:intrinsic-heat}
For every $w>0$,
\begin{equation}\label{Bcap}
 B(w)\le c(w)-\sqrt{\frac{x}{4\pi w}}\,(1-e^{-pw}).
\end{equation}
\end{lemma}
\begin{proof}
Write $k_s^{\rm round}$ and $k_s^M$ for the heat kernels with respect
to Riemannian volume on the sphere and on $M$. Fix $a\in M$ and put
\[
 v_s=\sqrt{4\pi s}\,k_s^{\rm round}(a,\cdot)\big|_M .
\]
For any smooth function $h$ on the sphere, minimality gives
\[
 \bar\Delta h=\Delta_M(h|_M)+\bar\nabla^2h(\nu,\nu)
 \qquad\text{on }M.
\]
Apply this identity to the spherical heat kernel. Hamilton's bound for
$\bar\nabla^2\log k_s^{\rm round}$ therefore gives
\[
 (\partial_s-\Delta_M)v_s
 =\sqrt{4\pi s}\left(\bar\nabla^2_{\nu,\nu} k_s^{\rm round}
                      +\frac{k_s^{\rm round}}{2s}\right)\ge0 .
\]

The measures $v_s\,d\sigma$ converge weakly to the unit point mass at $a$.
Indeed, in normal coordinates $y=\exp_a^M(u)$ the embedding is isometric
to first order at $a$, so
\[
 d_{\Sph}(a,y)^2=|u|^2+O(|u|^3),\qquad
 d\sigma(y)=(1+O(|u|^2))\,du,
\]
where $d_{\Sph}$ is the round distance. Set $u=\sqrt s\,\zeta$ in the
local heat kernel expansion. The limiting density is
\[
 \frac{1}{(4\pi)^{\frac n2}}\exp\!\left(-\frac{|\zeta|^2}{4}\right)\,d\zeta .
\]
In a sufficiently small coordinate ball the squared round distance is
at least $\tfrac12|u|^2$, so the Gaussian upper bound gives an
integrable majorant in $\zeta$. Outside
the coordinate neighborhood, compactness and embeddedness give a
positive round distance from $a$, so the integral there tends to zero
exponentially fast \cite{MP,NSS}. For $0<\varepsilon<s$, the heat equation
comparison principle gives
\[
 v_s\ge e^{(s-\varepsilon)\Delta_M}v_\varepsilon.
\]
Letting $\varepsilon\downarrow0$ gives $v_s\ge k_s^M(a,\cdot)$.

Now take $s=\tfrac{w}{n}$ and let $K_s^M$ be the heat kernel relative to $d\mu$.
After converting both kernels to the measures $d\pi$ and $d\mu$,
the bound $v_s\ge k_s^M(a,\cdot)$ becomes
\[
 K_s(a,b)\ge \sqrt{\frac{x}{4\pi w}}\,K_s^M(a,b)
 \qquad(a,b\in M).
\]
By symmetry,
\[
 B(w)=\int_M f^2 H_s\mu\,d\mu
 -\frac12\iint(f(a)-f(b))^2K_s(a,b)\,d\mu(a)d\mu(b).
\]
The first integral is at most $c(w)$ because $\int_M f^2\,d\mu=1$.
In the double integral the squared difference is nonnegative.
The kernel comparison therefore gives
\[
 \begin{aligned}
 \frac12\iint(f(a)-f(b))^2K_s(a,b)\,d\mu(a)d\mu(b)
 &\ge \sqrt{\frac{x}{4\pi w}}\,\langle f,(I-e^{s\Delta_M})f\rangle_{L^2(d\mu)}\\
 &=\sqrt{\frac{x}{4\pi w}}\,(1-e^{-pw}),
 \end{aligned}
\]
where the last equality uses $e^{s\Delta_M}f=e^{-pw}f$ and the normalization of $f$.
This proves \eqref{Bcap}.
\end{proof}

\section{Bounds for spherical harmonics of degrees two and three}\label{sec:coupled-moments}
We now estimate $a_2$ and $a_3$ in \eqref{green}. The estimate for $a_2$
uses the normal vector $\nu$. The estimate for $a_3$ also uses the
position vector $X$. To bound the two terms together, we use an identity
relating the averages of $XX^T$ and $\nu\nu^T$.

Let $\mathbb E$ denote integration over $M$ with respect to the probability
measure $\frac{d\mu}{\mu_0}$. Normalize the eigenfunction in this measure
by setting $e=\sqrt{\mu_0}f$. Then
\[
 \mathbb E e=0,\qquad \mathbb E e^2=1,\qquad \mathbb E(eX)=0.
\]
The matrix $P_\nu=\nu\nu^T$
is the orthogonal projection onto the normal line spanned by $\nu$.
Let $A_{\rm ext}$ agree with $A$ on $TM$ and vanish on
$\operatorname{span}\{X,\nu\}$. Write $Q_A$ for the normalized square
of this operator and $W_0$ for the average normal projection, so that
\[
 Q_A=\frac{A_{\rm ext}^2}{S},\qquad W_0=\mathbb E P_\nu.
\]
For $0<\eta<1$, define the matrix weight $W$ and define its squared
norm bound $c_\eta$ by
\[
                     W=(1-\eta)P_\nu+\eta Q_A,\qquad c_\eta=(1-\eta)^2+\eta^2(1-\rho) .
\]
The matrix $W$ is positive semidefinite. Since $A$ is trace-free,
\[
 \max_i\frac{\kappa_i^2}{S}\le 1-\rho.
\]
The eigenvalues of $Q_A$ are nonnegative
and sum to one, so $\tr Q_A^2\le1-\rho$. Also $P_\nu:Q_A=0$, where the
colon denotes the Euclidean inner product of tensors.
It follows that
\[
 |W|^2=(1-\eta)^2+\eta^2\tr(Q_A^2)\le c_\eta.
\]

The Gauss, Weingarten and Codazzi equations give
\[
 L\nu=S\nu,\qquad L P_\nu=2S(P_\nu-Q_A).
\]
Apply the product rule to $XX^T$. Since $LX=nX$ and
$I-XX^T-P_\nu$ is the projection onto $TM$, we get
\[
 L(XX^T)=2dXX^T+2P_\nu-2I.
\]
Since $M$ is closed, $\mathbb E(LP_\nu)=0$ and $\mathbb E(L(XX^T))=0$.
The first equality gives $\mathbb E Q_A=\mathbb E P_\nu=W_0$. The second
determines $\mathbb E(XX^T)$. In particular,
\[
 \mathbb E W=W_0,\qquad d\,\mathbb E(XX^T)=I-W_0,\qquad d=n+1 .
\]
We will estimate the two harmonic terms together. Let $K_{23}$ be
the set of nonnegative pairs $(r_2,r_3)$ that satisfy
\[
                  (1+z)^2r_2+r_3\le1+z^2\qquad(z\ge0).
\]
This set is closed, and the choice $z=0$ gives $r_2+r_3\le1$.
Hence $K_{23}$ is compact.

\begin{lemma}\label{lem:covariance}
Define the weighted moments $W_1,W_2$ and the residual $W_\perp$ by
\[
 W_1=\E(eW),\qquad W_2=\E(e^2W),\qquad W_\perp=eW-W_1-eW_2.
\]
Then the pair
\[
 (r_2,r_3)=\left(\frac{|W_1|^2}{c_\eta},\frac{\E|W_\perp|^2}{c_\eta}\right)
\]
belongs to $K_{23}$. Let $M_2$ and $M_3$ be the tensor moments
\[
 M_2=\E(eP_\nu),\qquad M_3=\E(eX\otimes P_\nu).
\]
With $\delta_3=1-\frac{\eta(n+\lambda)}{2S}$, we have
\begin{equation}\label{cov:main}
       d \delta_3^2|\operatorname{Sym}M_3|^2
       \le c_\eta r_3+\frac{\eta^2d\lambda}{2S^2}.
\end{equation}
Here $\operatorname{Sym}$ denotes orthogonal projection onto tensors
symmetric in all indices.
\end{lemma}
\begin{proof}
Take two independent points $a,b$ with law $\frac{d\mu}{\mu_0}$. Since
$W(a)$ and $W(b)$ are positive semidefinite, $W(a):W(b)\ge0$. Use this
quantity as a weight on the product measure. The integrals of $1$,
$e(a)e(b)$ and $e(a)^2e(b)^2$ with this weight are $|W_0|^2$, $|W_1|^2$
and $|W_2|^2$, respectively. Cauchy--Schwarz gives
$|W_1|^4\le |W_0|^2|W_2|^2$.

The functions $1$ and $e$ are orthonormal in $L^2\!\left(M,\frac{d\mu}{\mu_0}\right)$.
Project each entry of $W$ and $eW$ onto their span. The projected
matrices are $W_0+eW_1$ for $W$ and $W_1+eW_2$ for $eW$. Thus $W_\perp$ is the
part of $eW$ orthogonal to $1$ and $e$. Using $|W|^2\le c_\eta$ and $\mathbb E e^2=1$, we obtain
\[
 |W_0|^2+|W_1|^2\le c_\eta,\qquad
 |W_1|^2+|W_2|^2+\E|W_\perp|^2\le c_\eta .
\]
Multiply the first inequality by $z^2$ and add the second. For $z\ge0$,
\[
 2z|W_1|^2\le2z|W_0||W_2|\le z^2|W_0|^2+|W_2|^2.
\]
Therefore
\[
 (1+z)^2|W_1|^2+\E|W_\perp|^2\le(1+z^2)c_\eta.
\]
This proves that
$(r_2,r_3)\in K_{23}$.

To prove \eqref{cov:main}, write
\[
 W=P_\nu-\frac{\eta L P_\nu}{2S},\qquad \delta_2=1-\frac{\eta\lambda}{2S}.
\]
Integration by parts transfers $L$ from $P_\nu$ to $e$ or to $eX$.
Since $Le=\lambda e$ and $L(eX)=(n+\lambda)eX-2\nabla e$, we get
\[
 W_1=\delta_2M_2,\qquad
 \mathbb E(eX\otimes W)=\delta_3M_3+
       \frac{\eta}{S}\mathbb E(\nabla e\otimes P_\nu).
\]
Take a symmetric tensor $T$ of order three with $|T|=1$. Write
$T(v)=T(v,\cdot,\cdot)$ and $t_T=\mathbb E|T(\nu)|^2$.
The identity $d\,\mathbb E(XX^T)=I-W_0$ gives
$d\,\mathbb E|T(X)|^2=1-t_T$, so $0\le t_T\le1$.
Integrating $LX=nX$ gives $\mathbb E X=0$. Since also $\mathbb E(eX)=0$,
the terms $W_1$ and $eW_2$ both
have zero inner product with $T(X)$ after integration. Thus
\[
 |\mathbb E(eT(X):W)|=|\mathbb E(T(X):W_\perp)|
 \le\sqrt{\frac{(1-t_T)\mathbb E|W_\perp|^2}{d}}.
\]
For the term containing $\nabla e$, symmetry of $T$ gives
\[
 T(\nabla e,\nu,\nu)=T(\nu):\operatorname{Sym}(\nabla e\otimes\nu).
\]
Since $\nabla e\perp\nu$, we have
\[
 |\operatorname{Sym}(\nabla e\otimes\nu)|^2=\frac{|\nabla e|^2}{2}.
\]
Cauchy--Schwarz and $\mathbb E|\nabla e|^2=\lambda$ give
\[
 |\mathbb E T(\nabla e,\nu,\nu)|\le\sqrt{\frac{t_T\lambda}{2}}.
\]
Pair the identity for $\mathbb E(eX\otimes W)$ with $T$ and use the
two bounds. This gives
\[
 \sqrt d\,|\delta_3\,T:M_3|
 \le\sqrt{(1-t_T)\mathbb E|W_\perp|^2}
       +\frac{\eta}{S}\sqrt{\frac{t_Td\lambda}{2}}.
\]
The vector $(\sqrt{1-t_T},\sqrt{t_T})$ has norm one. Apply
Cauchy--Schwarz to the right side, square, and take the supremum over
all symmetric $T$ of norm one. This proves \eqref{cov:main}.
\end{proof}
For the rest of the proof we take
\[
 \eta=\frac13,\qquad c_\eta=\frac{5-\rho}{9}.
\]

\subsection{Estimating the two harmonic terms}
For $k=2,3$, each spherical harmonic of degree $k$ is the restriction
of a polynomial $P_T(X)=T(X,\ldots,X)$. Here $T$ is a constant
symmetric trace-free tensor of order $k$ on $\mathbb R^{d+1}$. Its
associated polynomial has squared norm
\[
 \|P_T\|_{L^2(\pi)}^2
 =\frac{k!\,|T|^2}{(d+1)(d+3)\cdots(d+2k-1)} .
\]
To see this, evaluate $P_T^2$ on a standard Gaussian vector. The Gaussian moment formula
sums over pairings of the $2k$ indices in the two copies of $T$.
A pairing contributes zero if it pairs two indices from the same copy,
because $T$ is trace-free. Each of the $k!$ pairings between the copies
contributes $|T|^2$. The
direction of the Gaussian vector is uniform on the sphere and is
independent of its length. Dividing by the $2k$-th moment of its
length gives the formula above.

The projection onto $TM$ is $I-XX^{\mathsf T}-P_\nu$. Take the trace of
the Euclidean Hessian of $P_T$ over $TM$ and use $\Delta_M X=-nX$.
Then
\[
 L P_T(X)=k(n+k-1)P_T(X)
            +k(k-1)T(X^{k-2},\nu,\nu).
\]
Here $X^{k-2}$ stands for $k-2$ copies of $X$. Multiply by $e$ and
integrate. Since $Le=\lambda e$, the identities for $k=2,3$ are
\[
 (2d-\lambda)\E(eP_T)=-2T:M_2,\qquad
 \bigl(3(d+1)-\lambda\bigr)\E(eP_T)=-6T:M_3 .
\]
Taking the supremum over all nonzero symmetric trace-free tensors $T$
of order $k$ gives
\[
 a_k=\sup_T\frac{\mu_0|\E(eP_T)|^2}{\|P_T\|_2^2}.
\]
Write $N_{d,k}=(d+1)(d+3)\cdots(d+2k-1)$ for the dimension-dependent
factor in the norm formula. For $k=2$, use
$W_1=\delta_2M_2$. For $k=3$, use $T:M_3=T:\operatorname{Sym}M_3$.
The norm formula for $P_T$ now gives
\[
 \frac{a_2}{\mu_0}\le\frac{2N_{d,2}|W_1|^2}{(2d-\lambda)^2\delta_2^2},
 \qquad
 \frac{a_3}{\mu_0}\le
 \frac{6N_{d,3}|\operatorname{Sym}M_3|^2}{\bigl(3(d+1)-\lambda\bigr)^2}.
\]
Multiply each bound by its weight in \eqref{green},
\[
 e^{-t\sigma_k}\left(\frac{1}{\sigma_k}-\frac{1}{\sigma_4}\right).
\]
Use $|W_1|^2=c_\eta r_2$ in the first bound and
\eqref{cov:main} in the second. After division by $\mu_0$, the sum of the degree-two and degree-three terms
in \eqref{green} is at most
\begin{equation}\label{eq:low}
       c_\eta(2\beta_2 r_2+\beta_3 r_3)+\beta_3\frac{\eta^2(1+\rho)p x^2}{2}.
\end{equation}
The coefficients are
\begin{align}
 \beta_2&=\frac{(1+6\rho)e^{-2(1+2\rho)t}}{4(2+2\rho-p)^2\delta_2^2},\nonumber\\
 \beta_3&=\frac{(1+7\rho)(1+2\rho)(1+6\rho)e^{-3(1+3\rho)t}}
 {2(1+\rho)(1+3\rho)(3+6\rho-p)^2\delta_3^2}.
 \label{eq:AB}
\end{align}
The same region $K_{23}$ applies to every choice of $n,S$ and $p$.

\subsection{Reduction to the coordinate eigenvalue}\label{sec:spectral-endpoint}
Set $t=\frac{x}{6}$ and fix $x$ and $\rho$. Let $U_{\rm flux}(p)$ be
the upper bound for $E_f^{-1}$ from the flux estimate \eqref{eq:naturalflux}.
To define $U_{\rm spec}(p)$, bound $B(w)$ in \eqref{green} by
\eqref{heatcap} and \eqref{Bcap}, including at $w=t$. For the last sum
in \eqref{green}, multiply the right side of \eqref{eq:low} by $c_I(D)$,
using $\mu_0\le c_I(D)$, and take the maximum over $(r_2,r_3)\in K_{23}$.
The sum of these bounds defines $U_{\rm spec}(p)$ for every real
$\tfrac12<p\le1$. When $p=\frac{\lambda}{n}$ comes from the eigenfunction
under consideration, $E_f\le U_{\rm spec}(p)$.
We will prove
\[
 U_{\rm flux}(p)U_{\rm spec}(p)\le U_{\rm flux}(1)U_{\rm spec}(1)\qquad\left(\frac12<p\le1\right).
\]

Set $\gamma=\frac{3x}{2}$. Differentiation gives
\[
 \frac{d}{dp}\log U_{\rm flux}(p)=
 \frac{2\bigl(1+\gamma(1-p-\varepsilon)\bigr)}
 {(2p-1+\varepsilon)(1+\gamma p)}>\tfrac12>t.
\]
The lower bound uses
\[
 p\le1,\qquad \gamma\le\frac32,\qquad
 \varepsilon\le\frac{\rho}{2}\le\frac18.
\]
Thus both $U_{\rm flux}(p)$ and $U_{\rm flux}(p)e^{-pw}$
increase with $p$ for $0<w\le t$.

Write $h_{\rm G}(w)$ for the Gaussian factor and $h_I(w)$ for the
bound from the normal collar. The upper bound for $B(w)$ is then
\[
 h_I(w)-h_{\rm G}(w)+h_{\rm G}(w)e^{-pw},
\]
where
\[
 h_{\rm G}(w)=\sqrt{\frac{x}{4\pi w}},\qquad
 h_I(w)=c_I\!\left(\left(\frac12+\rho\right)x+\frac{x}{4w}\right).
\]
Since
\[
 I\!\left(D+\frac{x}{4w}\right)
 \le\int_0^\infty \exp\!\left(-\frac{xz^2}{4w}\right)\,dz
 =\sqrt{\frac{\pi w}{x}},
\]
we have $h_I(w)\ge h_{\rm G}(w)$.
Multiplication by $U_{\rm flux}(p)$ therefore gives a sum of two nondecreasing
functions of $p$. Here $\delta_2=1-pt$ and $\delta_3=1-(1+p)t\ge1-2t>0$.
Thus all denominator factors in \eqref{eq:AB} that depend on $p$ are
positive and decrease with $p$. The coefficients in \eqref{eq:low} therefore
increase with $p$. The integral in $w$ and the maximum over $K_{23}$ both preserve
monotonicity in $p$. Hence
$U_{\rm flux}(p)U_{\rm spec}(p)\le U_{\rm flux}(1)U_{\rm spec}(1)$.

\section{Comparison at the coordinate eigenvalue}\label{sec:green-corners}
It remains to compare the two energy bounds at $p=1$, where $p=\tfrac{\lambda}{n}$. Recall that
\[
 x=\frac{n}{S},\qquad \rho=\frac{1}{n},\qquad x_*=\frac{7}{10}.
\]
Set
\[
 v=\frac{4\rho}{1+4\rho},\qquad s=\frac{1-v}{4}=\frac{1}{\sigma_4},
 \qquad t=\frac{x}{6}.
\]
The parameters lie in the rectangle $0\le x\le x_*$, $0\le v\le\tfrac12$.
At $x=0$ we use the limiting values of the estimates.
Write $\Gamma(x,v)=U_{\rm spec}(1)$ for the upper bound on $E_f$ at $p=1$.
For the reciprocal energy $E_f^{-1}$ the flux bound is
\[
 U_{\rm flux}(1)=\frac{1+\varepsilon}{\ell(x)}.
\]
Here $\ell$ is the affine function with $\ell(0)=\tfrac{4}{13}$ and
$\ell(x_*)=\tfrac{6}{25}(1+\tfrac{3x_*}{2})$. Thus it suffices to prove
$(1+\varepsilon)\Gamma(x,v)<\ell(x)$.

\paragraph{The heat integral.}
Recall that $c_I(y)=\tfrac{1}{2I(y)}$, where
$I(y)=\int_0^1 e^{-yz^2}\sqrt{1-z^2}\,dz$.
Set
\[
 h=\sqrt{\frac{3}{2\pi}},\qquad c_*=\frac{11}{4\pi}.
\]
Concavity of $c_I$ gives
\[
 c_I\!\left(\frac32+D\right)
 \le c_I(0)+\left(\frac32+D\right)c_I'(0)
 =c_*+\frac{D}{2\pi}\qquad(D\ge0).
\]
By Lemma~\ref{lem:heat-cap}, $(c_I^2)'\le\tfrac{1}{\pi}$. For
$0<z\le1$ and $D\ge0$, integration from $\tfrac32$ to
$\tfrac{3}{2z^2}+D$ gives
\[
 z^2c_I\!\left(\frac{3}{2z^2}+D\right)^2
 \le h^2+\left(c_*^2-h^2+\frac{D}{\pi}\right)z^2.
\]
Set $r=z^2$ in the integral below and apply Cauchy--Schwarz in $z$.
For the second inequality use $2\sqrt{y}\le\tfrac34+\tfrac{4y}{3}$.
\[
 \begin{aligned}
 \int_0^1 c_I\!\left(\frac{3}{2r}+D\right)\,dr
 &\le 2\sqrt{\frac{2h^2+c_*^2+\tfrac{D}{\pi}}{3}}\\
 &\le \frac34+\frac49\left(c_*^2+\frac{3+D}{\pi}\right).
 \end{aligned}
\]

\paragraph{An upper bound on the rectangle.}
Let $D=(\tfrac12+\rho)x$ as in Section~\ref{sec:heats}. The heat terms
in \eqref{green} are $\int_0^t B(w)\,dw+sB(t)$.
Set $w=tr$ in the integral. Then
\[
 \frac{x}{4w}=\frac{3}{2r},\qquad
 \sqrt{\frac{x}{4\pi w}}=\frac{h}{\sqrt r}.
\]
Apply \eqref{heatcap} and \eqref{Bcap}, and use
$1-e^{-tr}\ge tr(1-\tfrac{t}{2})$ for $0\le r\le1$.
The heat terms are at most
\[
 t\int_0^1 c_I\!\left(\frac{3}{2r}+D\right)\,dr
 +s\,c_I\!\left(\frac32+D\right)
 -h\left(1-\frac{t}{2}\right)\left(st+\frac{2t^2}{3}\right).
\]
The factor $\tfrac23$ comes from $\int_0^1\sqrt r\,dr$.

Put $\overline D=\tfrac{x(1+v)}{2}$. From $\rho\le\tfrac{v}{2}$ and
$s\rho=\tfrac{v}{16}$ we have
\[
 D\le \overline D,\qquad sD=x\left(\frac{s}{2}+\frac{v}{16}\right).
\]
The estimates for the integral and for $c_I(\tfrac32+D)$ now bound the heat terms by
\[
 \begin{aligned}
 U_{\rm heat}(x,v)={}&t\left(\frac34+\frac49\left(c_*^2+\frac{3+\overline D}{\pi}\right)\right)
 +sc_*+\frac{x}{2\pi}\left(\frac{s}{2}+\frac{v}{16}\right)\\
 &-h\left(1-\frac{t}{2}\right)\left(st+\frac{2t^2}{3}\right).
 \end{aligned}
\]
The tangent line to $c_I$ at zero also gives
\[
 \mu_0\le c_I(D)\le \overline\mu(x,v):=\frac{2}{\pi}\left(1+\frac{\overline D}{4}\right).
\]

It remains to bound the contributions from spherical harmonics of degrees
two and three.
For the coefficients of $r_2$ and $r_3$, respectively, set
\[
 \alpha_2(x,v)=\frac{(1+\tfrac{v}{2})e^{-2t}}{2(1-t)^2},\qquad
 \alpha_3(x,v)=\frac{(1+\tfrac{3v}{2})e^{-3t}}{8(1-2t)^2}.
\]
At $p=1$, we claim that $2\beta_2\le \alpha_2$ and $\beta_3\le \alpha_3$
in \eqref{eq:low}.
In \eqref{eq:AB}, use $\delta_2=1-t$, $\delta_3=1-2t$ and
$e^{-k(1+k\rho)t}\le e^{-kt}$ for $k=2,3$. The remaining factors satisfy
\[
 \begin{aligned}
 \frac{1+6\rho}{(1+2\rho)^2}&\le1+\frac{v}{2},\\
 \frac{(1+7\rho)(1+2\rho)(1+6\rho)}{(1+\rho)(1+3\rho)^3}
 &\le1+\frac{3v}{2}.
 \end{aligned}
\]
In each inequality, substitute $v=\tfrac{4\rho}{1+4\rho}$ and clear denominators.
After this multiplication, subtract the left side from the right side.
Every coefficient of the resulting polynomial in $\rho$ is nonnegative.
With $\eta=\tfrac13$ the bound \eqref{eq:low} has
$c_\eta=\tfrac{5-\rho}{9}\le\tfrac59$. Using $\rho\le\tfrac{v}{2}$, we obtain
\[
 \Gamma(x,v)\le\sup_{(r_2,r_3)\in K_{23}}\Psi(x,v;r_2,r_3),
\]
where
\[
 \Psi(x,v;r_2,r_3)=U_{\rm heat}+\overline\mu\left(\frac59(\alpha_2r_2+\alpha_3r_3)
       +\frac{x^2 \alpha_3}{18}\left(1+\frac{v}{2}\right)\right).
\]

\paragraph{Reduction to the four corners.}
Fix $(r_2,r_3)\in K_{23}$. We prove convexity of $\Psi$ in $v$ for
fixed $x$, and in $x$ for fixed $v$.
For fixed $x$, $U_{\rm heat}$ is affine in $v$ and the remaining term is a
polynomial in $v$ with nonnegative coefficients.
For fixed $v$, use $t=\tfrac{x}{6}$. As a polynomial in $t$, $U_{\rm heat}$ has a positive
cubic coefficient. Its quadratic coefficient is
\[
 \frac{4(1+v)}{3\pi}-\frac{13+3v}{24}h
 \ge (1+v)\left(\frac{4}{3\pi}-\frac{13h}{24}\right)>0.
\]
Thus $U_{\rm heat}$ is convex in $t$. The $t$-dependent factors in $\alpha_2,\alpha_3$ have the form
\[
 \frac{e^{-kt}}{(1-rt)^2},\qquad (k,r)=(2,1)\text{ or }(3,2).
\]
Their logarithmic derivatives $\tfrac{2r}{1-rt}-k$ are nonnegative and increasing.
The functions $\alpha_2$, $\alpha_3$, $\overline\mu$ and $x^2=36t^2$ are all
nonnegative, increasing and convex in $t$. Their products are therefore
convex in $t$. It follows that $\Psi$ is convex in $x=6t$.

\paragraph{The comparison at the corners.}
At each of the four corners, put $a=\alpha_2(x,v)$ and $b=\alpha_3(x,v)$.
Their formulas give $0<b\le a$. By the definition of $K_{23}$,
every pair $(r_2,r_3)\in K_{23}$ satisfies
\[
 (1+z)^2r_2+r_3\le 1+z^2\qquad(z\ge0).
\]
Take $z=\sqrt{\tfrac{a}{b}}-1$ and multiply by $b$. This gives
\[
 ar_2+br_3\le a+2b-2\sqrt{ab}.
\]
Thus $\Psi\le\Psi_{\rm c}(x,v)$ at each corner, where
\[
 \Psi_{\rm c}(x,v)=U_{\rm heat}
 +\overline\mu\left(
 \frac59\left(a+2b-2\sqrt{ab}\right)
 +\frac{x^2 b}{18}\left(1+\frac{v}{2}\right)\right).
\]
Use $t=0$ at $x=0$ and $t=\tfrac{7}{60}$ at $x=x_*$. Substitution gives
\[
 \max\!\left\{\Psi_{\rm c}(x,0),\frac98 \Psi_{\rm c}\!\left(x,\frac12\right)\right\}<\ell(x)
 \qquad(x=0,x_*).
\]

Let $\delta>0$ be the smallest of the four corner values of
$(1-\tfrac{2v}{9})\ell(x)-\Psi_{\rm c}(x,v)$. For each fixed pair
$(r_2,r_3)\in K_{23}$, convexity of $\Psi$ gives
\[
 \Psi(x,v;r_2,r_3)\le(1-\tfrac{2v}{9})\ell(x)-\delta
\]
throughout the rectangle. Indeed, the right side is affine in each
variable. First interpolate in $v$ on the edges $x=0,x_*$. Then
interpolate in $x$ for each fixed $v$. Since $\delta$ is independent of
$(r_2,r_3)$, the same bound holds for $\Gamma$.
Finally, the function
\[
 \frac{1}{1+\tfrac{\rho}{2}}=\frac{1-v}{1-\tfrac{7v}{8}}
\]
is concave in $v$. The line joining its values at $v=0$ and
$v=\tfrac12$ is $1-\tfrac{2v}{9}$. Therefore
\begin{equation}\label{elegant:strictmargin}
                   \Gamma(x,v)<\frac{\ell(x)}{1+\tfrac{\rho}{2}}.
\end{equation}

\begin{proof}[Proof of Theorem~\ref{main}]
Section~\ref{sec:reduction} proves the theorem for $n\le3$ and for $S\le n$.
In the remaining case, \eqref{eq:domain} holds. If the conclusion fails,
the same section gives $Lf=\lambda f$ with $0<\lambda\le n$, $f\perp V_X$
and $\int_M f^2\,d\mu=1$.
Section~\ref{sec:spectral-endpoint} gives $U_{\rm flux}(p)U_{\rm spec}(p)\le U_{\rm flux}(1)U_{\rm spec}(1)$.
Since $U_{\rm spec}(1)=\Gamma(x,v)$ and $\varepsilon\le\tfrac{\rho}{2}$,
\eqref{elegant:strictmargin} gives
\[
 1=E_f^{-1}E_f
 \le U_{\rm flux}(p)U_{\rm spec}(p)\le U_{\rm flux}(1)U_{\rm spec}(1)
 =\frac{(1+\varepsilon)\Gamma(x,v)}{\ell(x)}<1.
\]
Hence no eigenfunction with $0<\lambda\le n$ is orthogonal to $V_X$.
By Section~\ref{sec:reduction}, $\lambda_1=n$ and $\ker(L-n)=V_X$.
The dimension of $V_X$ gives the stated multiplicity.
\end{proof}

\section{A spectral estimate in terms of the principal curvatures}
\label{sec:principal-spectrum}

Let $M^n\subset\Sph^{n+1}(1)$ be a smooth closed connected embedded minimal
hypersurface, with $n\ge2$. Let $k$ be the maximum of the squared principal
curvatures on $M$, so that
\[
 k=\max_{p\in M}\|A_p\|_{\rm op}^2.
\]
In this section, orthogonality is taken in $L^2(M,d\sigma)$.
When $S$ is constant, Theorem~\ref{main} gives $\lambda_1=n$ and
identifies the first eigenspace with $V_X$. The following estimate then
applies to every eigenvalue above $n$.

\begin{proposition}\label{prop:principal-spectrum}
If $\lambda_1=n$, every eigenfunction $f\perp V_X$ with positive
eigenvalue $\lambda$ satisfies
\begin{equation}\label{eq:principal-spectrum}
 \lambda\ge 2(n+1-k).
\end{equation}
If $k\le\tfrac n2$, then $\lambda_1=n$ and $\ker(L-n)=V_X$, even when $S$
varies. In this case the smallest eigenvalue strictly greater than $n$
is at least $2(n+1-k)\ge n+2$.
\end{proposition}

\begin{proof}
We first allow any value of $\lambda_1$. Let $f$ be a real eigenfunction
with $Lf=\lambda f$, $\lambda>0$ and $f\perp V_X$. Normalize it by
$\int_M f^2\,d\sigma=1$. We view $\nabla f$ and $fX$ as
$\mathbb R^{n+2}$-valued functions on $M$. Both have mean zero. For $fX$
this follows from $f\perp V_X$. The $i$-th component of $\nabla f$ is
$\langle\nabla f,\nabla X^i\rangle$, and integration by parts gives
\[
 \int_M\langle\nabla f,\nabla X^i\rangle\,d\sigma
 =\int_M fLX^i\,d\sigma=n\int_M fX^i\,d\sigma=0.
\]

Write $D$ for Euclidean differentiation. For vector-valued functions
$U,V$ of mean zero, set
\[
 Q(U,V)=\int_M\bigl(\langle DU,DV\rangle-\lambda_1\langle U,V\rangle\bigr)\,d\sigma.
\]
Rayleigh's inequality applied to each component gives $Q(U,U)\ge0$.
For each tangent vector $v$, the Gauss formula gives
\[
 \begin{aligned}
 D_v(fX)&=v(f)X+fv,\\
 D_v\nabla f&=\nabla_v\nabla f+\langle Av,\nabla f\rangle\nu-v(f)X.
 \end{aligned}
\]
The Gauss equation gives $\operatorname{Ric}_M(v,v)=(n-1)|v|^2-|Av|^2$.
Since $\int_M|\nabla f|^2\,d\sigma=\lambda$, the integrated Bochner
formula \cite[Theorem~1]{Reilly} gives
\[
 \int_M|\nabla^2 f|^2\,d\sigma
 =\lambda(\lambda-n+1)+\int_M|A\nabla f|^2\,d\sigma.
\]
Taking inner products in the two derivative formulas and integrating, we obtain
\[
 \begin{aligned}
 Q(fX,fX)&=n+\lambda-\lambda_1,\\
 Q(\nabla f,fX)&=\int_M(f\Delta_M f-|\nabla f|^2)\,d\sigma=-2\lambda,\\
 Q(\nabla f,\nabla f)&=\lambda(\lambda-n+2-\lambda_1)
                           +2\int_M|A\nabla f|^2\,d\sigma.
 \end{aligned}
\]
Now $Q(fX,fX)>0$ since $\lambda\ge\lambda_1$. The Cauchy--Schwarz
inequality for $Q$ and the definition of $k$ give
\begin{equation}\label{eq:principal-rayleigh}
 n-2+\lambda_1-\lambda+\frac{4\lambda}{n+\lambda-\lambda_1}
 \le \frac{2}{\lambda}\int_M|A\nabla f|^2\,d\sigma\le 2k.
\end{equation}
Setting $\lambda_1=n$ gives \eqref{eq:principal-spectrum}.

Suppose now that $k\le\tfrac n2$. Since the coordinate functions have
eigenvalue $n$, we have $\lambda_1\le n$. If $\lambda_1<n$, every first
eigenfunction is orthogonal to $V_X$. Take $\lambda=\lambda_1$ in
\eqref{eq:principal-rayleigh}. By the strict bound
$\lambda_1>\tfrac n2$ \cite[Theorem~1.1]{DSS}, we obtain
\[
 2k\ge n-2+\frac{4\lambda_1}{n}>n.
\]
This contradicts $k\le\tfrac n2$. Hence $\lambda_1=n$.
Since $2(n+1-k)\ge n+2$, \eqref{eq:principal-spectrum} excludes any
nonzero eigenfunction at $n$ orthogonal to $V_X$. Thus $\ker(L-n)=V_X$.
Every eigenfunction with eigenvalue greater than $n$ is orthogonal to
$V_X$, so \eqref{eq:principal-spectrum} also gives the stated lower bound.
\end{proof}

The estimate \eqref{eq:principal-spectrum} is attained on the balanced
minimal Clifford product. Let $n=2m$ and take
\[
 M=\Sph^m\!\left(\frac{1}{\sqrt2}\right)
       \times\Sph^m\!\left(\frac{1}{\sqrt2}\right).
\]
Its principal curvatures are $1$ and $-1$, so $k=1$. Each factor has
first positive eigenvalue $n$. The product of a coordinate function from
each factor therefore has eigenvalue $2n$ and is orthogonal to $V_X$.
Thus it attains the lower bound $2(n+1-k)=2n$.


\begin{thebibliography}{99}

\bibitem{ABR}
S. Axler, P. Bourdon, and W. Ramey,
\emph{Harmonic function theory}, 2nd ed.,
Graduate Texts in Mathematics, vol.~137, Springer-Verlag, New York, 2001.
\doi{10.1007/978-1-4757-8137-3}.

\bibitem{Cecil}
T. E. Cecil, \emph{Isoparametric and Dupin hypersurfaces},
SIGMA Symmetry Integrability Geom. Methods Appl. \textbf{4} (2008),
Paper No.~062, 28 pp. \doi{10.3842/SIGMA.2008.062}.

\bibitem{Chang}
S. Chang, \emph{A closed hypersurface with constant scalar and mean
curvatures in $S^4$ is isoparametric}, Comm. Anal. Geom.
\textbf{1} (1993), no.~1, 71--100. \doi{10.4310/CAG.1993.v1.n1.a4}.

\bibitem{Chen}
R. Chen, \emph{An improved second pinching theorem for minimal
hypersurfaces with constant scalar curvature in spheres},
preprint, 2026. \arxiv{2609.17586v1}.

\bibitem{CLWY}
Q.-M. Cheng, F. Li, G. Wei, and J. Yin,
\emph{An improved result on Chern conjecture},
preprint, 2026. \arxiv{2610.05276v2}.

\bibitem{Chern}
S.-S. Chern,
\href{https://books.google.com/books?id=0u1NAQAAMAAJ}{\emph{Minimal
submanifolds in a Riemannian manifold}}, Technical Report No.~19,
Department of Mathematics, University of Kansas, Lawrence, KS, 1968,
58 pp.

\bibitem{CDK}
S.-S. Chern, M. do Carmo, and S. Kobayashi,
\emph{Minimal submanifolds of a sphere with second fundamental form
of constant length}, Functional Analysis and Related Fields
(Proc. Conf. in Honor of M.~H.~Stone, Chicago, IL, 1968)
(F.~E.~Browder, ed.), Springer-Verlag, Berlin, 1970, 59--75.
\doi{10.1007/978-3-642-48272-4_2}.

\bibitem{ChoeSoret}
J. Choe and M. Soret,
\emph{First eigenvalue of symmetric minimal surfaces in $\mathbb S^3$},
Indiana Univ. Math. J. \textbf{58} (2009), no.~1, 269--282.
\doi{10.1512/iumj.2009.58.3192}.

\bibitem{ChoiWang}
H. I. Choi and A.-N. Wang,
\emph{A first eigenvalue estimate for minimal hypersurfaces},
J. Differential Geom. \textbf{18} (1983), no.~3, 559--562.
\doi{10.4310/jdg/1214437788}.

\bibitem{DSS}
J. A. J. Duncan, Y. Sire, and J. Spruck,
\emph{An improved eigenvalue estimate for embedded minimal
hypersurfaces in the sphere},
Int. Math. Res. Not. IMRN (2024), no.~18, 12556--12567.
\doi{10.1093/imrn/rnae154}.

\bibitem{Hamilton}
R. S. Hamilton, \emph{A matrix Harnack estimate for the heat equation},
Comm. Anal. Geom. \textbf{1} (1993), no.~1, 113--126.
\doi{10.4310/CAG.1993.v1.n1.a6}.

\bibitem{Howard}
R. Howard, \emph{Blaschke's rolling theorem for manifolds with boundary},
Manuscripta Math. \textbf{99} (1999), no.~4, 471--483.
\doi{10.1007/s002290050186}.

\bibitem{IK}
A. Izosimov and B. Khesin, \emph{Vortex sheets and diffeomorphism groupoids},
Adv. Math. \textbf{338} (2018), 447--501. \doi{10.1016/j.aim.2018.09.015}.

\bibitem{JTZ}
A. Jim\'enez, C. Tapia Chinchay, and D. Zhou,
\emph{A lower bound for the first eigenvalue of a minimal hypersurface in
the sphere}, Rev. Mat. Iberoam. \textbf{42} (2026), no.~1, 261--278.
\doi{10.4171/RMI/1587}.

\bibitem{KM}
R. Kusner and P. McGrath, \emph{On Steklov eigenspaces for free boundary
minimal surfaces in the unit ball}, Amer. J. Math. \textbf{146} (2024),
no.~5, 1275--1293. \doi{10.1353/ajm.2024.a937942}.

\bibitem{LXX}
L. Lei, H. Xu, and Z. Xu, \emph{On Chern's conjecture for minimal hypersurfaces in spheres},
preprint, 2017. \arxiv{1712.01175}.

\bibitem{Lawson}
H. B. Lawson, Jr., \emph{Local rigidity theorems for minimal hypersurfaces},
Ann. of Math. (2) \textbf{89} (1969), no.~1, 187--197.
\doi{10.2307/1970816}.

\bibitem{MP}
S. Minakshisundaram and {\AA}. Pleijel,
\emph{Some properties of the eigenfunctions of the Laplace-operator on
Riemannian manifolds}, Canad. J. Math. \textbf{1} (1949), no.~3,
242--256. \doi{10.4153/CJM-1949-021-5}.

\bibitem{MOU}
H. Muto, Y. Ohnita, and H. Urakawa,
\emph{Homogeneous minimal hypersurfaces in the unit spheres and the first eigenvalues of their Laplacian},
Tohoku Math. J. (2) \textbf{36} (1984), no.~2, 253--267. \doi{10.2748/tmj/1178228851}.

\bibitem{NSS}
A. Nowak, P. Sj\"ogren, and T. Z. Szarek,
\emph{Sharp estimates of the spherical heat kernel},
J. Math. Pures Appl. (9) \textbf{129} (2019), 23--33.
\doi{10.1016/j.matpur.2018.10.002}.

\bibitem{PengTerng}
C.-K. Peng and C.-L. Terng,
\emph{Minimal hypersurfaces of spheres with constant scalar curvature},
Seminar on Minimal Submanifolds (E.~Bombieri, ed.), Ann. of Math. Stud.,
vol.~103, Princeton University Press, Princeton, NJ, 1983, 177--198.
\doi{10.1515/9781400881437-010}.

\bibitem{Reilly}
R. C. Reilly, \emph{Applications of the Hessian operator in a Riemannian
manifold}, Indiana Univ. Math. J. \textbf{26} (1977), no.~3, 459--472.
\doi{10.1512/iumj.1977.26.26036}.

\bibitem{Simons}
J. Simons, \emph{Minimal varieties in Riemannian manifolds},
Ann. of Math. (2) \textbf{88} (1968), no.~1, 62--105.
\doi{10.2307/1970556}.

\bibitem{SuhYang}
Y. J. Suh and H. Y. Yang, \emph{The scalar curvature of minimal hypersurfaces in a unit sphere},
Commun. Contemp. Math. \textbf{9} (2007), no.~2, 183--200. \doi{10.1142/S021919970700237X}.

\bibitem{Takahashi}
T. Takahashi, \emph{Minimal immersions of Riemannian manifolds},
J. Math. Soc. Japan \textbf{18} (1966), no.~4, 380--385.
\doi{10.2969/jmsj/01840380}.

\bibitem{TangYan}
Z. Tang and W. Yan, \emph{Isoparametric foliation and Yau conjecture
on the first eigenvalue}, J. Differential Geom. \textbf{94} (2013),
no.~3, 521--540. \doi{10.4310/jdg/1370979337}.

\bibitem{Yau}
S.-T. Yau, \emph{Problem section}, Seminar on Differential Geometry
(S.-T.~Yau, ed.), Ann. of Math. Stud., vol.~102,
Princeton University Press, Princeton, NJ, 1982, 669--706.
\doi{10.1515/9781400881918-035}.

\bibitem{Yu}
J. Yu, \emph{On the first eigenvalue of embedded minimal hypersurfaces
in the unit sphere}, preprint, 2026. \arxiv{2606.10962v2}.

\bibitem{Zeng}
L. Zeng, \emph{The first eigenvalue of embedded minimal hypersurfaces
in the unit sphere I: Yau's conjecture}, preprint, 2025,
34 pp. \arxiv{2508.06123v1}.

\bibitem{Zhao}
Y. Zhao, \emph{The first eigenvalue of embedded minimal hypersurfaces
in the unit sphere}, preprint, 2026. \arxiv{2603.20890v2}.

\end{thebibliography}
\end{document}